\documentclass[11pt,oneside,reqno]{amsart}
\usepackage[all]{xy}
\usepackage{amsfonts,amsmath,amssymb,amscd}
\usepackage[breaklinks]{hyperref}
\usepackage{caption} 
\usepackage{mathrsfs}
\usepackage[top=1in, bottom=1in, left=1in, right=1in]{geometry}

\renewcommand{\a}{\alpha}
\renewcommand{\b}{\beta}

 \renewcommand{\a}{\alpha}
\renewcommand{\b}{\beta}

\newcommand{\g}{\gamma}
\newcommand{\G}{\Gamma}
\renewcommand{\l}{\lambda}

\renewcommand{\(}{\left\(}
\renewcommand{\)}{\right\)}
\renewcommand{\[}{\left\[}
\renewcommand{\]}{\right\]}
\renewcommand{\i}{\infty}

\numberwithin{equation}{section}

\newtheorem{theorem}{Theorem}

\newtheorem{remark}[subsection]{Remark}
\newtheorem{remark*}[subsection]{Remark}

\newtheorem{corollary}[theorem]{Corollary}

\numberwithin{theorem}{section}

   \makeatletter
\def\proof{\@ifnextchar[{\@oproof}{\@nproof}}
\def\@oproof[#1][#2]{\trivlist\item[\hskip\labelsep\textit{#2 Proof of\
#1.}~]\ignorespaces}
\def\@nproof{\trivlist\item[\hskip\labelsep\textit{Proof.}~]\ignorespaces}

\makeatother

\begin{document}
\title{Extended higher Herglotz functions I.	Functional equations}

\author{Atul Dixit, Rajat Gupta and Rahul Kumar}
\address{Discipline of Mathematics, Indian Institute of Technology, Gandhinagar, Palaj, Gandhinagar 382355, Gujarat, India}\email{adixit@iitgn.ac.in; rajat\_gupta@iitgn.ac.in; rahul.kumr@iitgn.ac.in}


\maketitle
\pagenumbering{arabic}
\pagestyle{headings}
\begin{center}
\dedicatory{\emph{Dedicated to Professor Don Zagier on the occasion of his 70th birthday}}
\end{center}
\begin{abstract}
In 1975, Don Zagier obtained a new version of the Kronecker limit formula for a real quadratic field which involved an interesting function $F(x)$ which is now known as the \emph{Herglotz function}. As demonstrated by Zagier, and very recently by Radchenko and Zagier, $F(x)$ satisfies beautiful properties which are of interest in both algebraic number theory as well as in analytic number theory. In this paper, we study $\mathscr{F}_{k,N}(x)$, an extension of the Herglotz function which also subsumes \emph{higher Herglotz function} of Vlasenko and Zagier. We call it the \emph{extended higher Herglotz function}. It is intimately connected with a certain generalized Lambert series. We derive two different kinds of functional equations satisfied by $\mathscr{F}_{k,N}(x)$. Radchenko and Zagier gave a beautiful relation between the integral $\displaystyle\int_{0}^{1}\frac{\log(1+t^x)}{1+t}\, dt$ and $F(x)$ and used it to evaluate this integral at various rational as well as irrational arguments. We obtain a relation between $\mathscr{F}_{k,N}(x)$ and a generalization of the above integral involving polylogarithm. The asymptotic expansions of $\mathscr{F}_{k, N}(x)$ and some generalized Lambert series are also obtained along with other supplementary results.

\end{abstract}{}

\tableofcontents

\section{Introduction}\label{intro}

The Kronecker limit formula is one of the celebrated theorems in number theory having applications not only in other areas of Mathematics such as in geometry but also in Physics, for example, in string theory for the one-loop computation in 
the Polyakov's perturbative approach \cite{takhtajan}. Let $K$ denote a number field and $\zeta_{K}(s)$ the Dedekind zeta-function of $K$. If $A$ denotes an ideal class group of $K$, then we can write $\zeta_{K}(s)=\sum_{A}\zeta(s, A),$ where, for Re$(s)>1$, $\zeta(s, A)=\sum_{a\in A}\frac{1}{N(a)^s},$ with $N(a)$ being the norm of the fractional ideal $a$. Then the Kronecker limit formula is concerned with the evaluation of the constant term $\varrho(A)$ in the Laurent series expansion of $\zeta(s, A)$. This was done by Kronecker \cite{kronecker} in the case of an imaginary quadratic field $K$. While Hecke \cite{hecke} obtained a corresponding formula for a real quadratic field $K$, it was not explicit as one of the expressions in it was left as an integral containing the Dedekind eta-function. 

In a seminal paper \cite{zagier1975}, Zagier obtained an explicit form of the Kronecker limit formula for real quadratic fields. Zagier's formula gives a new representation for the value
\begin{equation*}
\varrho(\mathcal{B}):=\lim_{s\to1}\left(D^{s/2}\zeta(s, \mathcal{B})-\frac{\log\epsilon}{s-1}\right),
\end{equation*}
where $\mathcal{B}$ is an element of the group of narrow ideal classes of invertible fractional $\mathcal{O}_{D}$-ideals\footnote{Here two ideals belong to the same narrow class if their quotient is a principal ideal $\lambda\mathcal{O}_{D}$ with $N(\l)>0$.} of the quadratic order $\mathcal{O}_{D}=\mathbb{Z}+\mathbb{Z}\frac{D+\sqrt{D}}{2}$ of discriminant $D>0$, $\zeta(s, \mathcal{B})$ is the sum of $N(\mathfrak{a})^{-s}$ over all invertible $\mathcal{O}_{D}$-ideals $\mathfrak{a}$ in the class $\mathcal{B}$, and $\epsilon=\epsilon_{D}$ is the smallest unit $>1$ in $\mathcal{O}_{D}$ of norm $1$. The reason such a representation is worth studying is because the special value at $s=1$ of the Dirichlet series $L_{K}(s, \chi)=\sum_{\mathfrak{a}}\chi(\mathfrak{a})N(\mathfrak{a})^{-s}$ equals $D^{-1/2}\sum_{\mathcal{B}}\chi(\mathcal{B})\varrho(\mathcal{B})$ for any character $\chi$ on the narrow ideal class group. Zagier's formula \cite[p.~164]{zagier1975} reads
\begin{equation}\label{pww}
\varrho(\mathcal{B})=\sum_{w\in\text{Red}(\mathcal{B})}P(w,w'),
\end{equation}
where the function $P(x, y)$ for $x>y>0$ is defined as
\begin{equation}\label{pxy}
P(x, y):=F(x)-F(y)+\textup{Li}_{2}\left(\frac{y}{x}\right)-\frac{\pi^2}{6}+\log\left(\frac{x}{y}\right)\left(\gamma-\frac{1}{2}\log(x-y)+\frac{1}{4}\log\left(\frac{x}{y}\right)\right),
\end{equation}
with $\text{Red}(\mathcal{B})$ being the set of larger roots $w=(-b+\sqrt{D})/(2a)$ of all reduced primitive quadratic forms $Q(X, Y)=aX^2+bXY+cY^2 (a, c>0, a+b+c<0)$ of discriminant $D$ belonging to the class $\mathcal{B}$. Here $\textup{Li}_{2}(t)$ is the dilogarithm function 
\begin{align}\label{li2t}
\sum_{n=1}^{\infty}\frac{t^n}{n^2}\hspace{5mm}(0<t<1),
\end{align}
and 
\begin{equation}\label{herglotzdef}
F(x):=\sum_{n=1}^{\infty}\frac{\psi(nx)-\log(nx)}{n}\hspace{5mm}\left(x\in\mathbb{C}\backslash(-\infty,0]\right),
\end{equation} 
where $\psi(a)=\G'(a)/\G(a)$ is the logarithmic derivative of the gamma function. The function $F(x)$ is similar to a function that Herglotz studied in \cite{herglotz} and hence Radchenko and Zagier \cite{raza} call it the \textit{Herglotz function}.

The Herglotz function $F(x)$ satisfies several interesting and important properties as has been demonstrated by Zagier in \cite{zagier1975}, and very recently by Radchenko and Zagier in \cite{raza} where they study, among other things, the relation of this function with the Dedekind eta-function, functional equations satisfied by $F(x)$ in connection with Hecke operators, the cohomological aspects of $F(x)$, and its special values at positive rationals and quadratic units. Prior to this, Zagier \cite[Equations (7.4), (7.8)]{zagier1975} had obtained beautiful functional equations that $F(x)$ satisfies, namely, for $x\in\mathbb{C}\backslash(-\infty,0]$,
{\allowdisplaybreaks\begin{align}
F(x)-F(x+1)-F\left(\frac{x}{x+1}\right)&=-F(1)+\textup{Li}_2\left(\frac{1}{1+x}\right),\label{fe1}\\
F(x)+F\left(\frac{1}{x}\right)&=2F(1)+\frac{1}{2}\log^{2}(x)-\frac{\pi^2}{6x}(x-1)^2,\label{fe2}
\end{align}}
where \cite[Equation (7.12)]{zagier1975}
\begin{equation*}
F(1)=-\frac{1}{2}\g^2-\frac{\pi^2}{12}-\g_1,
\end{equation*}
and $\g$ and $\g_1$ denote Euler's constant and the first Stieltjes constant respectively. 

Apart from their intrinsic beauty, these functional equations are very useful, for example, Zagier \cite{zagier1975} used them to prove Meyer's theorem and Radchenko and Zagier \cite{raza} have used them to obtain asymptotic expansions of $F(x)$ as $x\to0$ and $x\to 1$. Functional equations of the type in \eqref{fe1} occur in a variety of areas of mathematics such as period functions for Maass forms \cite[Equation (0.1)]{lewzag} (see also \cite[Equation (1)]{betcon}), cotangent functions, double zeta functions etc. We refer the reader to a beautiful article of Zagier \cite{zagier3term} on this topic.

Several generalizations of the Herglotz function (also known sometimes as the \emph{Herglotz-Zagier function}) have been studied in the literature. Ishibashi \cite{ishibashi} studied the $j^{\textup{th}}$ order Herglotz function defined by
\begin{equation*}
\Phi_{j}(x):=\sum_{n-1}^{\infty}\frac{R_j'(nx)-\log^{j}(nx)}{n} \hspace{7mm}(j\in\mathbb{N}),
\end{equation*}
where $R_j(x)$ is a generalization of the function $R(x)$ studied by Deninger \cite{deninger}, and is defined as the solution to the difference equation
\begin{equation*}
f(x+1)-f(x)=\log^{j}(x),\hspace{5mm}(f(1)=(-1)^{j+1}\zeta^{(j)}(0)),
\end{equation*}
with $\zeta(s)$ denoting the Riemann zeta function. It is clear that $\Phi_1(x)=F(x)$. Ishibashi's function $\Phi_{j}(x)$ appears in his explicit representation \cite[Theorem 3]{ishibashi} of the Laurent series coefficients of $Z_{Q_{\ell}}(s)$, the zeta function associated to indefinite quadratic forms $Q_{\ell}$, where
\begin{equation*}
D^{s/2}\zeta(s, \mathcal{B})=\sum_{\ell}Z_{Q_{\ell}}(s)\hspace{5mm} (\textup{Re}(s)>1).
\end{equation*}
Zagier \cite{zagier1975} had earlier obtained the constant coefficient and the coefficient of $(s-1)^{-1}$ in the Laurent series expansion of $D^{s/2}\zeta(s, \mathcal{B})$. The Herglotz function $F(x)$ makes its conspicuous presence in Zagier's expression of the constant coefficient as can be seen from \eqref{pww} and \eqref{pxy}.

Masri \cite{masri} generalized $F(x)$ in a different direction by considering
\begin{equation*}
F(\a, s, x):=\sum_{n=1}^{\infty}	\frac{\a(n)\left(\psi(nx)-\log(nx)\right)}{n^s}\hspace{5mm}(\textup{Re}(s)>0, x>0),
\end{equation*}
where $\a:\mathbb{Z}\to\mathbb{C}$ is a periodic function with period $M$. Vlasenko and Zagier \cite{vz} studied the following special case of Masri's $F(\a, s, x)$ when $\a(n)\equiv1$ and $s>1$ is a natural number but without considering $\log(nx)$ term (which is not needed for convergence of the series since $s>1$): 
\begin{equation}\label{vlaszag0}
F_k(x):=\sum_{n=1}^{\infty}\frac{\psi(nx)}{n^{k}}\hspace{5mm}\left(k\in\mathbb{N},k>1, x\in\mathbb{C}\backslash(-\infty,0]\right).
\end{equation}
They termed it the \emph{higher Herglotz function}\footnote{It is important to note that Vlasenko and Zagier \cite{vz} define their higher Herglotz function by $\sum_{n=1}^{\infty}\psi(nx)/n^{k-1}$ for $k>2$.}
and obtained following beautiful functional relations for it \cite[Equations (11), (12)]{vz}:
{\allowdisplaybreaks\begin{align}\label{vz1f}
F_k(x)+(-x)^{k-1}F_k\left(\frac{1}{x}\right)&=-\gamma\zeta(k)\left(1+(-x)^{k-1}\right)-\sum_{r=2}^{k-1}\zeta(r)\zeta(k+1-r)(-x)^{r-1}\nonumber\\
&\quad+\zeta(k+1)\left((-x)^{k}-\frac{1}{x}\right),
\end{align}}
and\footnote{The interpretation of $\zeta(1,k-1)$ given in \cite[p.~28]{vz} has a minus sign missing in front of the whole term, that is, the correct form is $-\left(\zeta(k-1,1)+\zeta(k)-\g\zeta(k-1)\right)$.} 
\begin{align}\label{vz2f}
F_k(x)-F_k(x+1)+(-x)^{k-1}F_k\left(\frac{x+1}{x}\right)&=(-x)^{k-1}\left(\zeta(k,1)+\zeta(k+1)-\g\zeta(k)\right)\nonumber\\
&\quad-\sum_{r=1}^{k-1}\zeta(k+1-r,r)(-x)^{r-1}+\zeta(k+1)\left(\frac{(-x)^{k}}{x+1}-\frac{1}{x}\right),
\end{align}
where
\begin{equation*}
\zeta(m, n):=\sum_{p>q>0}\frac{1}{p^mq^n}\hspace{5mm}(m\geq2, n\geq1).
\end{equation*}
The function $F_k(x)$ plays an instrumental role in the higher Kronecker ``limit'' formulas developed in \cite{vz}.

The goal of our paper is to study a novel generalization of $F(x)$, which we call the \emph{extended higher Herglotz function}. We define it by
\begin{align}\label{fknx}
\mathscr{F}_{k,N}(x):=\sum_{n=1}^{\infty}\frac{\psi(n^{N}x)-\log(n^{N}x)}{n^k}\hspace{5mm}\left(x\in\mathbb{C}\backslash(-\infty,0]\right),
\end{align}
where $k$ and $N$ are positive real numbers such that $k+N>1$. It is clear that when $k>1$, one can write
\begin{align}\label{fknxk>1}
\mathscr{F}_{k,N}(x)&=\sum_{n=1}^{\infty}\frac{\psi(n^{N}x)}{n^{k}}-\sum_{n=1}^{\infty}\frac{\log(n^{N}x)}{n^k}\nonumber\\
&=\sum_{n=1}^{\infty}\frac{\psi(n^{N}x)}{n^{k}}+N\zeta'(k)-\zeta(k)\log x.
\end{align}
However, we keep \eqref{fknx} as the definition of the extended higher Herglotz function so as to keep many special cases under one roof. Indeed, from \eqref{vlaszag0} and \eqref{fknxk>1}, we have for $k\in\mathbb{N}$,
\begin{equation}\label{vlaszag}
\mathscr{F}_{k,1}(x)=F_k(x)+\zeta'(k)-\zeta(k)\log x.
\end{equation}
Also $\mathscr{F}_{1,1}(x)$ is the Herglotz function defined in \eqref{herglotzdef}.

We were naturally led to the idea of considering the extended higher Herglotz function in \eqref{fknx} from our previous work \cite{dgkm} on a generalized Lambert series. Indeed, the extended higher Herglotz function turns up after differentiating a transformation for the generalized Lambert series 
\begin{equation}\label{gls}
\sum_{n=1}^{\infty}\frac{n^{-2Nm-1}\textup{exp}\left(-a(2n)^{N}\alpha\right)}{1-\textup{exp}\left(-(2n)^{N}\alpha\right)}
\end{equation}
with respect to $a$ and then letting $a=1$ as can be seen from \eqref{zetageneqna} below. 

An application of our extended higher Herglotz functions is in obtaining the asymptotic expansion of the Lambert series 
\begin{equation*}
\sum_{n=1}^{\infty}\frac{n^{-2Nm-1+N}}{\textup{exp}\left((2n)^{N}\alpha\right)-1}\hspace{5mm}(N\hspace{1mm}\text{odd})
\end{equation*}
as $\a\to0^{+}$ which can be seen in Theorem \ref{lambert asymptotic}. Obtaining asymptotic expansions for Lambert series and its generalizations has received a lot of attention since last century because of their applications in the theory of the Riemann zeta function. For example, Wigert \cite{wig0} (see also \cite[p.~163, Theorem 7.15]{titch}) obtained the asymptotic expansion of $\sum_{n=1}^{\infty}1/(e^{nz}-1)$ as $z\to0$ which was used by Lukkarinen \cite{mari} to not only study the moments of the Riemann zeta function but to also obtain meromorphic continuation of the modified Mellin transform of $\left|\zeta(1/2+ix)\right|^{2}$ defined for $\textup{Re}(s)>1$ by $\int_{1}^{\infty}\left|\zeta\left(\frac{1}{2}+ix\right)\right|^{2}x^{-s}\, dx$. For odd $N\in\mathbb{N}$, Wigert himself \cite[p.~9-10]{wigk} obtained the asymptotic expansion of the generalized Lambert series 
\begin{equation*}
\sum_{n=1}^{\infty}\frac{1}{e^{n^Nx}-1}=\sum_{j=1}^{\infty}\sum_{n=1}^{\infty}e^{-n^{N}jx}.
\end{equation*}
as $x\to0$. Note that the inner sum in the above double series is a generalized partial theta function. Wigert not only wrote a paper on this partial theta function \cite{Wig} but also on the above generalized Lambert series \cite{wigk}.

Zagier \cite{zagierasym} had conjectured that the asymptotic expansion of the Lambert series associated with the square of the Ramanujan tau function $\tau(n)$, that is of $\sum_{n=1}^{\infty}\tau^{2}(n)e^{-nx}$, as $x\to 0^{+}$, can be expressed in terms of the non-trivial zeros of the Riemann zeta function. This conjecture was proved by Hafner and Stopple \cite{hafstop}. Thus obtaining asymptotic expansions of Lambert series is a useful endeavor.

In this paper, we obtain two different kinds of functional equations satisfied by $\mathscr{F}_{k,N}(x)$. See Theorems \ref{hhfe} and \ref{genherglotzNodd}. As a special case of one of these functional equations, we get \eqref{fe2} of Zagier.

For $x>0$, one of the beautiful identities that $F(x)$ satisfies \cite{raza} is 
\begin{align}\label{zaginteval}
J(x)=F(2x)-2F(x)+F\left(\frac{x}{2}\right)+\frac{\pi^2}{12x},
\end{align}
where
\begin{equation}\label{jx}
J(x):=\int_0^1 \frac{\log\left(1+t^x\right)}{1+t}\ dt.
\end{equation}
Radchenko and Zagier \cite{raza} have obtained special values of $F(x)$ at positive rationals and quadratic units with the help of which they are able to evaluate special values of $J(x)$, for example, 
{\allowdisplaybreaks\begin{align}\label{jxspl}
J(4+\sqrt{17})&= -\frac{\pi^2}{6}+\frac{1}{2}\log^{2}(2)+\frac{1}{2}\log(2)\log(4+\sqrt{17}),\nonumber\\
J\left(\frac{2}{5}\right)&=\frac{11\pi^2}{240}+\frac{3}{4}\log^{2}(2)-2\log^{2}\left(\frac{\sqrt{5}+1}{2}\right).\nonumber
\end{align}} 
Earlier, Herglotz \cite{herglotz} as well as Muzaffar and Williams \cite{muzwil} had evaluated such integrals but only of a particular type, that is, $J(n+\sqrt{n^2-1})$. For example, Herglotz \cite[Equation (70a)]{herglotz} showed that
\begin{align*}
J(4+\sqrt{15})=-\frac{\pi^2}{12}(\sqrt{15}-2)+\log(2)\log(\sqrt{3}+\sqrt{5})+\log\left(\frac{\sqrt{5}+1}{2}\right)\log(2+\sqrt{3}).
\end{align*}
Regarding this integral, Chowla \cite[p.~372]{chowla} remarked, \textit{`A direct evaluation of this definite integral is probably difficult'}. Indeed, these evaluations come from not only employing methods from analytic number theory but also those from algebraic number theory, which is why they are all the more interesting!

In this paper, a generalization of \eqref{zaginteval} involving $\mathscr{F}_{k,N}(x)$ is obtained in Theorem \ref{zagintevalgen}. Lastly, we obtain some results showcasing the intimate relationship between generalized Lambert series and extended higher Herglotz functions.

%

\section{Main results}\label{mr}
Throughout this paper, the notation $\sum_{j=-(N-1)}^{(N-1)}{\vphantom{\sum}}''$ will denote a sum over over $j=-(N-1), -(N-3),\cdots, N-3, N-1$.

\subsection{Functional equations for extended higher Herglotz functions $\mathscr{F}_{k,N}(x)$}

\begin{theorem}\label{hhfe}
Let $k,N\in \mathbb{N}$ such that $1<k\leq N$ and $x\in\mathbb{C}\backslash(-\infty,0]$. Let
\begin{equation}\label{bknx}
\mathscr{B}(k, N, x):=\left\{
	\begin{array}{ll}
		(-1)^{k+N+1}x^{1/N}\zeta\left( 1+\frac{k}{N}\right) & \mbox{if } k=N \\
		0 & \mbox{if } k\neq N
	\end{array}
\right. .
\end{equation}
Then
\begin{align}\label{hhfeeqn}
&x^{\frac{1-k}{N}}\mathscr{F}_{k,N}(x)-\frac{(-1)^k}{N}\sum_{j=-(N-1)}^{(N-1)}{\vphantom{\sum}}''e^{i\pi j (k-1)/N}\mathscr{F}_{\frac{N+k-1}{N},\frac{1}{N}}\left(\frac{e^{-\frac{i\pi j}{N}}}{x^{1/N}}\right)\nonumber\\
&=\frac{\pi}{N} \frac{\zeta\left(1+\frac{k-1}{N}\right)}{\sin \left(\frac{\pi}{N}(k-1) \right)}+x^{\frac{1-k}{N}}\left(-\left(\gamma+\log x \right) \zeta(k)+N\zeta'(k)\right)-\frac{1}{x^{\frac{N+k-1}{N}}}\zeta(k+N)+\mathscr{B}(k, N, x).
\end{align}
Also,
\begin{align}\label{hhfeeqn1}
&\mathscr{F}_{1,N}(x)+\frac{1}{N}\sum_{j=-(N-1)}^{(N-1)}{\vphantom{\sum}}''\mathscr{F}_{1,\frac{1}{N}}\left(\frac{e^{-\frac{i\pi j}{N}}}{x^{1/N}}\right)\nonumber\\
&=\frac{1}{N}\left\{\frac{\pi^2}{6}-(N-1)\g\log x+\frac{1}{2}(\log x)^{2}-N\g^2-(N^2+1)\g_1\right\}-\frac{\zeta(N+1)}{x}+\mathscr{B}(1,N,x).
\end{align}
\end{theorem}
We note that when $N=1$, \eqref{hhfeeqn1} gives Zagier's \eqref{fe2} as a special case.

In the next theorem, we obtain a different kind of functional equation for $\mathscr{F}_{k,N}(x)$, more specifically, for a combination of $\mathscr{F}_{k,N}$ functions. There is a trade-off between the two types of functional equations obtained in Theorems \ref{hhfe} and \ref{genherglotzNodd}, namely, Theorem \ref{hhfe} is valid for any natural numbers $k$ and $N$ such that $1<k\leq N$ whereas Theorem \ref{genherglotzNodd} is valid for any \emph{odd} natural numbers $k$ and $N$. These two types are, indeed, equivalent when $k$ and $N$ are odd such that $1<k\leq N$. This is shown in subsection \ref{equivalence}. 

The idea for deriving the second type of functional equation, considered in Theorem \ref{genherglotzNodd} below, stemmed from a beautiful result obtained in \cite[Theorem 2.2]{dgkm}. This identity, which was instrumental in deriving the main transformation for a generalized Lambert series in \cite{dgkm}, is as follows.

Let $u\in\mathbb{C}$ be fixed such that \textup{Re}$(u)>0$. Then, 
\begin{equation}\label{raabesumeqn}
\sum_{m=1}^{\infty}\int_{0}^{\infty}\frac{t\cos(t)}{t^2+m^2u^2}\, \mathrm{d}t=\frac{1}{2}\left\{\log\left(\frac{u}{2\pi}\right)-\frac{1}{2}\left(\psi\left(\frac{iu}{2\pi}\right)+\psi\left(\frac{-iu}{2\pi}\right)\right)\right\}.
\end{equation}

\begin{theorem}\label{genherglotzNodd}
Let $k$ and $N$ be odd natural numbers. Let $\textup{Re}(x)>0$. Then,
\begin{align}\label{thmnotequaleqn}
&\sum_{j=-(N-1)}^{N-1}{\vphantom{\sum}}''\exp\left(-\tfrac{i\pi(k-1)j}{2N}\right)\left\{\mathscr{F}_{\frac{N+k-1}{N},\frac{1}{N}}\left(i\left(\tfrac{2\pi}{x}\right)^{\frac{1}{N}}e^{\frac{i\pi j}{2N}}\right)+\mathscr{F}_{\frac{N+k-1}{N},\frac{1}{N}}\left(-i\left(\tfrac{2\pi}{x}\right)^{\frac{1}{N}}e^{\frac{i\pi j}{2N}}\right)\right\}\nonumber\\
&=(-1)^{\frac{k+1}{2}}N\left(\tfrac{2\pi}{x}\right)^{\frac{k-1}{N}}\bigg\{\mathscr{F}_{k,N}\left(\tfrac{ix}{2\pi}\right)+\mathscr{F}_{k,N}\left(-\tfrac{ix}{2\pi}\right)+2\left(\left(\gamma-\log\left(\tfrac{2\pi}{x}\right)\right)\zeta(k)-N\zeta'(k)\right)+\mathscr{C}(k, N, x)\bigg\},
\end{align}
where, for $\frac{1-k}{N}\neq-2\left\lfloor\frac{k}{2N}\right\rfloor$,
\begin{align}\label{cknx1}
\mathscr{C}(k, N, x):=
\frac{\pi\zeta\left(\frac{N+k-1}{N}\right)}{N\sin\left(\frac{\pi}{2N}(1-k)\right)}\left(\frac{x}{2\pi}\right)^{\frac{k-1}{N}}+4\pi\sum_{j=1}^{\lfloor\frac{k}{2N}\rfloor}(-1)^j(2\pi)^{-2j-1}\zeta(k-2Nj)\zeta(2j+1)x^{2j},
\end{align}
and for $\frac{1-k}{N}=-2\left\lfloor\frac{k}{2N}\right\rfloor$,
\begin{align}\label{cknx2}
\mathscr{C}(k, N, x)&:=\frac{2(-1)^{\frac{k-1}{2N}}}{N}\left(\frac{x}{2\pi}\right)^{\frac{k-1}{N}}\left[\left(\gamma+\log\left(\tfrac{2\pi}{x}\right)\right)\zeta\left(1+\tfrac{k-1}{N}\right)-\zeta'\left(1+\tfrac{k-1}{N}\right)\right]\nonumber\\
&\quad+4\pi\sum_{j=1}^{\lfloor\frac{k}{2N}\rfloor-1}(-1)^j(2\pi)^{-2j-1}\zeta(k-2Nj)\zeta(2j+1)x^{2j}.
\end{align}
\end{theorem}
The remaining case not considered in the above theorem, namely, $\frac{1-k}{N}=-2\left\lfloor\frac{k}{2N}\right\rfloor=0$, that is, $k=1$, is given in the following result.
\begin{theorem}\label{thmtriplepole}
Let $k$ and $N$ be odd natural numbers and $\textup{Re}(x)>0$. Then,
\begin{align*}
&\sum_{j=-(N-1)}^{N-1}{\vphantom{\sum}}''e^{\frac{i\pi }{2}(j-1)}\left\{\mathscr{F}_{1,\frac{1}{N}}\left(i\left(\tfrac{2\pi}{x}\right)^{\frac{1}{N}}e^{\frac{i\pi j}{2N}}\right)+\mathscr{F}_{1,\frac{1}{N}}\left(-i\left(\tfrac{2\pi}{x}\right)^{\frac{1}{N}}e^{\frac{i\pi j}{2N}}\right)\right\}\nonumber\\
&=-N\left\{\mathscr{F}_{1,N}\left(\tfrac{ix}{2\pi}\right)+\mathscr{F}_{1,N}\left(-\tfrac{ix}{2\pi}\right)-\frac{1}{N}\left(\frac{\pi^2}{12}-2N\gamma^2+2(N-1)\gamma\log\left(\frac{2\pi}{x}\right)\right.\right.\nonumber\\
&\quad\left.\left.+\log^2(2\pi)-\log\left(\frac{4\pi^2}{x}\right)\log x-2(N^2+1)\gamma_1\right)\right\}.
\end{align*}
\end{theorem}

Ramanujan's famous formula for $\zeta(2m+1)$ \cite[p.~173, Ch. 14, Entry 21(i)]{ramnote}, \cite[pp.~319-320, formula (28)]{lnb}, \cite[pp.~275-276]{bcbramsecnote} is given for $\a, \b>0$, $\textup{Re}(\a)>0, \textup{Re}(\b)>0, \a\b=\pi^2$ and $m\in\mathbb{Z}\backslash\{0\}$, by
\begin{align}\label{zetaodd}
\a^{-m}\left\{\frac{1}{2}\zeta(2m+1)+\sum_{n=1}^{\infty}\frac{n^{-2m-1}}{e^{2\a n}-1}\right\}&=(-\b)^{-m}\left\{\frac{1}{2}\zeta(2m+1)+\sum_{n=1}^{\infty}\frac{n^{-2m-1}}{e^{2\b n}-1}\right\}\nonumber\\
&\quad-2^{2m}\sum_{j=0}^{m+1}\frac{(-1)^jB_{2j}B_{2m+2-2j}}{(2j)!(2m+2-2j)!}\a^{m+1-j}\b^j,
\end{align}
where $B_n$ is the $n^{\textup{th}}$ Bernoulli number. The Lambert series in \eqref{gls} enables us to obtain a two-parameter generalization of \eqref{zetaodd} given in \cite[Theorem 2.4]{dgkm}. See \eqref{zetageneqna} below.

As a special case, Theorem \ref{genherglotzNodd} gives a beautiful transformation for a combination of the Vlasenko-Zagier higher Herglotz function $F_k(x)$ which is analogous to \eqref{zetaodd}.
\begin{corollary}\label{trans2m+1}
Let $\alpha$ and $\beta$ be two complex numbers with $\textup{Re}(\a)>0, \textup{Re}(\b)>0$ and $\alpha\beta=4\pi^2$. Then for $m\in\mathbb{N}$,
\begin{align}\label{trans2m+1eqn}
&\alpha^{-m}\left\{2\gamma\zeta(2m+1)+\sum_{n=1}^\infty\frac{1}{n^{2m+1}}\left(\psi\left(\frac{in\alpha}{2\pi}\right)+\psi\left(-\frac{in\alpha}{2\pi}\right)\right)\right\}\nonumber\\
&=-(-\beta)^{-m}\left\{2\gamma\zeta(2m+1)+\sum_{n=1}^\infty\frac{1}{n^{2m+1}}\left(\psi\left(\frac{in\beta}{2\pi}\right)+\psi\left(-\frac{in\beta}{2\pi}\right)\right)\right\}\nonumber\\
&\qquad-2\sum_{j=1}^{m-1}(-1)^j\zeta(1-2j+2m)\zeta(2j+1)\alpha^{j-m}\beta^{-j}.
\end{align}
\end{corollary}
It is important to note here that Corollary \ref{trans2m+1} cannot be obtained for any general $\a, \b$ such that $\a\b=4\pi^2$ by applying Vlasenko and Zagier's \eqref{vz1f} once with $k=2m+2$ and $x=i\a/(2\pi)$, and again with $k=2m+2$ and $x=i\b/(2\pi)$, and adding the resulting two equations, for, doing so will result in
\begin{align}\label{zagiervla}
&\sum_{n=1}^\infty\frac{\psi\left(\frac{in\alpha}{2\pi}\right)}{n^{2m+1}}+(-\beta)^m\alpha^{-m}\sum_{n=1}^\infty\frac{\psi\left(-\frac{in\alpha}{2\pi}\right)}{n^{2m+1}}+\sum_{n=1}^\infty\frac{\psi\left(\frac{in\beta}{2\pi}\right)}{n^{2m+1}}+(-\alpha)^m\beta^{-m}\sum_{n=1}^\infty\frac{\psi\left(-\frac{in\beta}{2\pi}\right)}{n^{2m+1}}\nonumber\\
&=-\gamma\left(2+(-\beta)^m\alpha^{-m}+(-\alpha)^m\beta^{-m}\right)\zeta(2m+1)+i\zeta(2m+2)\Bigg(\frac{2\pi}{\beta}+\frac{2\pi}{\alpha}-\frac{(-\beta)^m\alpha^{-m}\beta}{2\pi}\nonumber\\
&\quad-\frac{(-\alpha)^m\beta^{-m}\alpha}{2\pi}\Bigg)-\sum_{r=2}^{2m}\zeta(r)\zeta(2m+2-r)\left(-\frac{i}{2\pi}\right)^{r-1}\left(\alpha^{r-1}+\beta^{r-1}\right).
\end{align}
It is clear that it is impossible to get Corollary \ref{trans2m+1} from \eqref{zagiervla} for any $\a,\b$ such that $\a\b=4\pi^2$ because of the factors of the form $(-\beta)^m\alpha^{-m}$ and $(-\alpha)^m\beta^{-m}$ in front of the series $\sum_{n=1}^\infty n^{-2m-1}\psi\left(-\frac{in\alpha}{2\pi}\right)$ and $\sum_{n=1}^\infty n^{-2m-1}\psi\left(-\frac{in\beta}{2\pi}\right)$ respectively in \eqref{zagiervla} whereas we do not have such factors in Corollary \ref{trans2m+1}. 

Only in the case when $\a=\b=2\pi$ and $m$ is replaced by $2m$ is Corollary \ref{trans2m+1} obtainable from \eqref{zagiervla}. We state below this special case of Corollary \ref{trans2m+1}, namely for $m\in\mathbb{N}$,
\begin{align}\label{trans4m+1eqn}
\sum_{n=1}^\infty\frac{1}{n^{4m+1}}\left(\psi(in)+\psi(-in)\right)&=-2\gamma\zeta(4m+1)-\sum_{j=1}^{2m-1}(-1)^{j}\zeta(2j+1)\zeta(4m+1-2j).
\end{align}

Letting $m=1$ in Corollary \ref{trans2m+1} gives a beautiful modular relation:
\begin{corollary}\label{modular}
Let $\alpha$ and $\beta$ be two complex numbers with $\textup{Re}(\a)>0, \textup{Re}(\b)>0$ and $\alpha\beta=4\pi^2$. Then
\begin{align}\label{modulareqn}
\frac{1}{\a}\left\{2\g\zeta(3)+\sum_{n=1}^{\infty}\frac{1}{n^{3}}\left(\psi\left(\frac{in\a}{2\pi}\right)+\psi\left(\frac{-in\a}{2\pi}\right)\right)\right\}=\frac{1}{\b}\left\{2\g\zeta(3)+\sum_{n=1}^{\infty}\frac{1}{n^{3}}\left(\psi\left(\frac{in\b}{2\pi}\right)+\psi\left(\frac{-in\a}{2\pi}\right)\right)\right\}.
\end{align}
\end{corollary}
\begin{remark}\label{ehhfgls1}
Corollary \textup{\ref{trans2m+1}} is an analogue of Ramanujan's formula for $\zeta(2m+1)$, that is, \eqref{zetaodd}. Indeed, the combination of Vlasenko-Zagier higher Herglotz functions in Corollary \textup{\ref{trans2m+1}} seems to mimick the role played by the Lambert series in Ramanujan's formula. 

This phenomenon is seen to be true even when we consider the generalization of Corollary \textup{\ref{trans2m+1}} given in Theorem \textup{\ref{genherglotzNodd}}. Indeed, the combination of the extended higher Herglotz functions $\mathscr{F}_{k, N}(x)$ in Theorem \textup{\ref{genherglotzNodd}} plays a similar role as that played by the generalized Lambert series in Theorem \textup{1.2} of \cite{dixitmaji1}. The latter theorem is the special case $a=1$ of \eqref{zetageneqna}.
\end{remark}
\subsection{Asymptotics of the extended higher Herglotz functions}\label{asymptotics}
An immediate application of the functional equations given in the above theorems is in obtaining the asymptotic expansions of the extended higher Herglotz functions $\mathscr{F}_{k,N}(x)$ as $x\to0$. These asymptotic expansions for different conditions on the parameters involved are collected together in the theorem below.
\begin{theorem}\label{asymptotic}
Let $x\in\mathbb{C}\backslash(-\infty,0]$ and $k, N$ are positive real numbers such that $k+N>1$. Let $\mathscr{B}(k, N, x)$ be defined in \eqref{bknx}.
Then as $x\to\infty$,
\begin{align}\label{fknxinfty}
\mathscr{F}_{k,N}(x) \sim-\frac{1}{2x}\zeta(k+N)-\sum_{n=1}^{\infty}\frac{B_{2n}}{2n}\zeta(k+2nN)x^{-2n}.
\end{align}
Also, as $x\to 0$,
{\allowdisplaybreaks\begin{align}
&\textup{(i)}\hspace{1mm}\text{For}\hspace{1mm} k,N\in\mathbb{N}\hspace{1mm}\text{such that}\hspace{1mm}1<k\leq N,\hspace{1mm}\text{we have}\nonumber\\
&\mathscr{F}_{k,N}(x) \sim -\frac{\zeta(k+N)}{x}-\left(\gamma+\log x \right) \zeta(k)+N\zeta'(k)+\frac{\pi}{N} \frac{\zeta\left(1+\frac{k-1}{N}\right)}{\sin \left(\frac{\pi}{N}(k-1) \right)}x^{\frac{k-1}{N}}+x^{\frac{k-1}{N}}\mathscr{B}(k, N, x)\nonumber\\
&\qquad\qquad+(-1)^k \sum_{m=1}^{\i}(-1)^{m(N+1)}\zeta(k-mN)\zeta\left(1+m\right)x^{m};\label{f1nx}\\
&\textup{(ii)}\hspace{1mm}\text{For}\hspace{1mm}N\in\mathbb{N},\nonumber\\
&\mathscr{F}_{1,N}(x)\sim-\frac{\zeta(N+1)}{x}+\frac{1}{N}\left\{\frac{\pi^2}{6}-(N-1)\g\log x+\frac{1}{2}\log^{2}x-N\g^2-(N^2+1)\g_1\right\}\nonumber\\
&\qquad\qquad\qquad+\mathscr{B}(1, N, x)-\sum_{m=1}^{\i}(-1)^{m(N+1)}\zeta(1-mN)\zeta\left(1+m\right)x^{m}\label{ff1nx}.
\end{align}}
\end{theorem}
Letting $N=1$ gives Radchenko and Zagier's asymptotic expansions for $F(x)$ as $x\to\infty$ and as $x\to0$ as can be seen from \cite[Equation (7)]{raza}.

\begin{theorem}\label{asymthm2.3}
Let $k$ and $N$ be positive real numbers such that $k+N>1$ and let \textup{Re}$(x)>0$. Let $\mathscr{C}(k, N, x)$ be defined in \eqref{cknx1} and \eqref{cknx2}. Then as $x\to\infty$,
\begin{align}\label{f+finf}
\mathscr{F}_{k, N}\left(\frac{ix}{2\pi}\right)+\mathscr{F}_{k, N}\left(-\frac{ix}{2\pi}\right)\sim\sum_{n=1}^{\infty}\frac{(-1)^{n+1}}{n}B_{2n}\zeta(k+2nN)\left(\frac{2\pi}{x}\right)^{2n},
\end{align}
Also, for odd natural numbers $k$ and $N$, as $x\to 0$,
\begin{align}\label{f+f}
&\mathscr{F}_{k, N}\left(\frac{ix}{2\pi}\right)+\mathscr{F}_{k, N}\left(-\frac{ix}{2\pi}\right)\nonumber\\
&\sim2\left\{\left(\log\left(\frac{2\pi}{x}\right)-\g\right)\zeta(k)+N\zeta'(k)\right\}-\mathscr{C}(k, N, x)+\frac{(-1)^{\frac{N+k+2}{2}}}{N}\left(\frac{x}{2\pi}\right)^{\frac{k-1}{N}}\exp{\left(-\frac{i\pi(N+k-1)}{2N}\right)}\nonumber\\
&\quad\times\sum_{n=1}^{\infty}\frac{B_{2n}}{n\cos\left(\frac{\pi n}{N}\right)}\zeta\left(\frac{N+k+2n-1}{N}\right)\left(\frac{x}{2\pi}\right)^{\frac{2n}{N}},
\end{align}
and for an odd natural number $N$, as $x\to 0$,
\begin{align}\label{f+f1}
&\mathscr{F}_{1, N}\left(\frac{ix}{2\pi}\right)+\mathscr{F}_{1, N}\left(-\frac{ix}{2\pi}\right)\nonumber\\
&\sim\frac{1}{N}\left(\frac{\pi^2}{12}-2N\gamma^2+2(N-1)\gamma\log\left(\frac{2\pi}{x}\right)+\log^2(2\pi)-\log\left(\frac{4\pi^2}{x}\right)\log x-2(N^2+1)\gamma_1\right)\nonumber\\
&\quad+\frac{i}{N}(-1)^{\frac{N+1}{2}}\sum_{n=1}^{\infty}\frac{B_{2n}}{n\cos\left(\frac{\pi n}{N}\right)}\zeta\left(1+\frac{2n}{N}\right)\left(\frac{x}{2\pi}\right)^{\frac{2n}{N}}.
\end{align}
\end{theorem}

\subsection{An interesting relation between a generalized polylogarithm function and the extended higher Herglotz functions}
The polylogarithm is a generalization of the dilogarithm function of \eqref{li2t}. It is defined by
\begin{align}
\textup{Li}_{s}(t):=\sum_{n=1}^{\infty}\frac{t^n}{n^s}.\nonumber
\end{align}
It is clear that $\textup{Li}_{s}(t)$ converges for any complex $s$ as long as $|t|<1$. It can be analytically continued for $|t|\geq1$. There are a number of generalizations of polylogarithm in the literature.

In the following result, we encountered a new generalization of the polylogarithm function. We define it for $N\in\mathbb{N}, s\in\mathbb{C}$ and $|t|<1$ by
\begin{equation}\label{lintN}
{}_N\textup{Li}_{s}(t):=\sum_{n=1}^{\infty}\frac{t^{n^{N}}}{n^s}.
\end{equation}
Clearly, ${}_1\textup{Li}_{s}(t)=\textup{Li}_{s}(t)$ for $|t|<1$.
\begin{theorem}\label{zagintevalgen}
Let $k$ and $N$ be any positive integers. Let ${}_N\textup{Li}_{s}(t)$ be defined in \eqref{lintN}. Define
\begin{align}\label{defnJkN}
J_{k,N}(x):=\int_0^1\left(\frac{2^{k-1}}{u-1}-\frac{2^Nu^{(2^N-1)}}{u^{2^N}-1}-\frac{2^{k-1}-1}{u\log u}\right){}_N\mathrm{Li}_k(-u^x) \ du.
\end{align}
Then for $x>0$, we have
\begin{align}\label{zagintevalgeneqn}
J_{k,N}(x)&=\mathscr{F}_{k,N}(2^Nx)-\left(2^{k-1}+2^{1-k}\right)\mathscr{F}_{k,N}(x)+\mathscr{F}_{k,N}\left(\frac{x}{2^N}\right)\nonumber\\
&\quad+\left(2^N+2^{-N}-2^{k-1}-2^{1-k}\right)\frac{\zeta(k+N)}{x}.
\end{align}
\end{theorem}
The integral in \eqref{defnJkN} is a natural generalization of \eqref{jx} in the setting of $\mathscr{F}_{k, N}(x)$ in the sense that it arises naturally while extending the Radchenko-Zagier relation between $J(x)$ and $F(x)$ in \eqref{zaginteval}. Indeed, it can be easily checked that $J_{1,1}(x)=J(x)$ and that \eqref{zaginteval} follows by letting $k=N=1$ in Theorem \ref{zagintevalgen}.

The above theorem leads to the following relation involving polylogarithm and Vlasenko and Zagier's higher Herglotz function $F_k(x)$ from \eqref{vlaszag}, which is valid for a positive integer $k>1$:
\begin{corollary}\label{N=1JK}
Let $x>0$ and let $k>1$ be a positive integer. Let $F_k(x)$ be defined in \eqref{vlaszag0}. Then
\begin{align}
&\int_0^1\left(\frac{2^{k-1}}{u-1}-\frac{2u}{u^{2}-1}-\frac{2^{k-1}-1}{u\log u}\right)\mathrm{Li}_k(-u^x) \ du\nonumber\\
&=F_k(2x)-\left(2^{k-1}+2^{1-k}\right)F_{k}(x)+F_{k}\left(\frac{x}{2}\right)+\left(2-2^{k-1}-2^{1-k}\right)\left(\zeta'(k)-\zeta(k)\log x+\frac{1}{x}\zeta(k+1)\right)\nonumber\\
&\quad+\frac{1}{2x}\zeta(k+1).\nonumber
\end{align}
\end{corollary}
For $k$ odd, Corollary \ref{N=1JK}, in turn, gives an evaluation of an integral involving polylogarithm function in terms of the higher Herglotz function $F_k(2)$, special values of the Riemann zeta function and its derivatives as stated below.
\begin{corollary}\label{N=1x1}
Let $k>1$ be an odd positive integer. Let $F_k(x)$ be defined in \eqref{vlaszag0}. Then
\begin{align}\label{N=1x1eqn}
&\int_0^1\left(\frac{2^{k-1}}{u-1}-\frac{2u}{u^{2}-1}-\frac{2^{k-1}-1}{u\log u}\right)\mathrm{Li}_k(-u) \ du\nonumber\\
&=\left(1-2^{1-k}\right)F_k(2)+(2^{k-1}-1)\g\zeta(k)+\left(\frac{1}{2}-2^{-k}\right)\zeta(k+1)+\left(2-2^{k-1}-2^{1-k}\right)\zeta'(k)\nonumber\\
&\quad+\sum_{r=2}^{k-1}(-1)^{r-1}\left(2^{k-2}+2^{-k}-2^{r-k}\right)\zeta(r)\zeta(k+1-r).
\end{align}
\end{corollary}

\subsection{Transformations involving higher Herglotz functions and generalized Lambert series}

In \cite[Corollary 2.13]{dkk}, the following companion of Ramanujan's formula \eqref{zetaodd} was obtained.
Let $m\in\mathbb{N}$. If $\a$ and $\b$ are complex numbers such that $\textup{Re}(\a)>0$, $\textup{Re}(\b)>0$, and $\a\b=\pi^2$, then
\begin{align}\label{companioneqn}
&\a^{-\left(m-\frac{1}{2}\right)}\left\{\frac{1}{2}\zeta(2m)+\sum_{n=1}^\infty\frac{n^{-2m}}{e^{2n\a}-1}\right\}-\sum_{j=0}^{m-1}\frac{2^{2j-1}B_{2j}}{(2j)!}\zeta(2m-2j+1)\a^{2j-m-\frac{1}{2}}\nonumber\\
&=(-1)^{m+1}\b^{-\left(m-\frac{1}{2}\right)}\left\{\frac{\g}{\pi}\zeta(2m)+\frac{1}{2\pi}\sum_{n=1}^\infty n^{-2m}\left(\psi\left(\frac{in\b}{\pi}\right)+\psi\left(-\frac{in\beta}{\pi}\right)\right)\right\}.
\end{align}

In what follows, we give a one-parameter extension of \eqref{companioneqn}. 
\begin{theorem}\label{gencompanion}
Let $N$ be an odd positive integer. For $\a \b^N= \pi^{N+1}$ with $\textup{Re}(\a)>0, ~\textup{Re}(\b) >0$, and $m \geq 1$,
\begin{align}\label{gencompanioneqn}
&\a^{-\left(\frac{2Nm}{N+1}-\frac{1}{2}\right)}\left(\frac{1}{2}\zeta(2Nm+1-N)+\sum_{n=1}^{\i}\frac{n^{-2Nm-1+N}}{e^{(2n)^N\a}-1} \right)\nonumber\\
&\hspace{2cm}-\sum_{j=0}^{m-1}\frac{B_{2j}~\zeta(2Nm+1-2Nj)}{(2j)!}2^{N(2j-1)}\a^{2j-\frac{2Nm}{N+1}-\frac{1}{2}}\nonumber\\
&=\frac{2^{2m(N-1)}}{N\pi^{\frac{N+1}{2}}}(-1)^{m+1}\b^{-\left(\frac{2Nm}{N+1}-\frac{N}{2}\right)}\Bigg[\frac{N\gamma}{2^{N-1}}\zeta(2m)\nonumber\\
&\hspace{2cm}+\frac{1}{2^N}\sum_{j=-\frac{(N-1)}{2}}^{\frac{(N-1)}{2}}\sum_{n=1}^{\i}\frac{1}{n^{2m}}\left(\psi\left( \frac{i \b}{2\pi}(2n)^{1/N}e^{\frac{i\pi j}{N}}\right)+\psi\left( -\frac{i \b}{2\pi}(2n)^{1/N}e^{\frac{i\pi j}{N}}\right) \right) \Bigg].
\end{align}
\end{theorem}
It is straightforward to see that letting $N=1$ in the above theorem gives \eqref{companioneqn}.

\begin{theorem}\label{gendivdkk}
Let $N$ be an odd positive integer. If $\a,\b>0$ such that $\a\b^{N}=\pi^{N+1}$,
\begin{align}\label{gendivdkkeqnalphabeta}
&\sum_{n=1}^\infty\frac{n^{N-1}}{\exp\left((2n)^N\alpha\right)-1}-\frac{1}{2^N\alpha N}\left(N\gamma-\log(2\pi)-(N-1)(\log2-\gamma)\right)\nonumber\\
&=\frac{1}{2^N\alpha(N+1)}\log\left(\frac{\alpha}{\beta}\right)+\frac{2}{2^N\alpha N}\sum_{j=-\left(\frac{N-1}{2}\right)}^{\frac{N-1}{2}}\sum_{n=1}^\infty\Bigg[\log\left(\frac{\beta}{2\pi}(2n)^{\frac{1}{N}}e^{\frac{i\pi j}{N}}\right)\nonumber\\
&\quad-\frac{1}{2}\left\{\psi\left(i\frac{\beta}{2\pi}(2n)^{\frac{1}{N}}e^{\frac{i\pi j}{N}}\right)+\psi\left(-i\frac{\beta}{2\pi}(2n)^{\frac{1}{N}}e^{\frac{i\pi j}{N}}\right)\right\}\Bigg]+\begin{cases}
\frac{\alpha}{2},&\ \mathrm{if},\ N=1\\
0,&\ \mathrm{if},\ N>1.
\end{cases}
\end{align}
%
\end{theorem}

As an application of the asymptotic expansions derived in Theorem \ref{asymthm2.3}, we obtain the complete asymptotic expansions for the Lambert series with the help of Theorem \ref{gencompanion} above.

\begin{theorem}\label{lambert asymptotic}
Let $N$ be an odd positive integer and let $m \geq 1$. As $\a \to 0$,
\begin{align}\label{lambert asymptotic eqn1}
\sum_{n=1}^{\i}&\frac{n^{-2Nm-1+N}}{e^{(2n)^N\a}-1}\sim \frac{1}{2^{N}\a}\zeta(2Nm+1)-\frac{1}{2}\zeta(2Nm+1-N)+\frac{(-1)^{m+1}}{N\pi^{2m}}2^{(2m-1)(N-1)}\a^{2m-1}\nonumber\\
&\qquad\qquad\qquad\times\sum_{\ell=1}^{\infty}\frac{(-1)^{\ell+1}B_{2\ell N}}{\ell}\zeta\left(2m+2\ell\right)\left( \frac{2^{N-1}\a}{\pi}\right)^{2\ell}+\mathscr{D}_{N}(m,\a),
\end{align}
where
\begin{align}
&\mathscr{D}_{N}(m,\a):=\frac{2^{(2m-1)(N-1)}}{N\pi^{2m}}(-1)^{m+1}\a^{2m-1}\Bigg\{\zeta(2m)\left(N\gamma-\log\left( \frac{\a2^{N-1}}{\pi}\right)\right)-\zeta'(2m)\Bigg\}\nonumber\\
&\qquad \qquad \qquad \quad+\sum_{j=1}^{m-1}\frac{B_{2j}~\zeta(2Nm+1-2Nj)}{(2j)!}2^{N(2j-1)}\a^{2j-1}.\nonumber
\end{align}

Also, in particular, as $\a \to 0$,
{\allowdisplaybreaks\begin{align}\label{lambert asymptotic eqn2}
&\sum_{n=1}^{\i}\frac{n^{-2m}}{e^{2n\a}-1}\sim \frac{\zeta(2m+1)}{2\a}-\frac{1}{2}\zeta(2m)+\frac{(-1)^{m+1}\a^{2m-1}}{\pi^{2m}}\Bigg\{\zeta(2m)\left(\gamma-\log\left( \frac{\a}{\pi}\right)\right)-\zeta'(2m)\Bigg\}\nonumber\\
&\qquad\qquad\qquad+\frac{(-1)^{m+1}\a^{2m-1}}{\pi^{2m}}\sum_{\ell=1}^{\infty}\frac{(-1)^{\ell+1}B_{2\ell}}{\ell}\zeta\left(2m+2\ell\right)\left( \frac{\a}{\pi}\right)^{2\ell}\nonumber\\
&\qquad\qquad\qquad+\sum_{j=1}^{m-1}\frac{B_{2j}~\zeta(2m+1-2j)}{(2j)!}(2\a)^{(2j-1)}.
\end{align}}
\end{theorem}

We note in passing that the well-known result \cite[p.~903, Formula \textbf{8.361.8}]{grad}, for $\textup{Re}(x)>0$, namely
\begin{align}
\psi(x)-\log x=-\int_0^\infty\left(\frac{1}{1-e^{-t}}-\frac{1}{t}\right)e^{-xt}\ dt,\nonumber
\end{align}
gives an integral representation for $\mathscr{F}_{k,N}(x)$:
\begin{align}\label{intreprF}
\mathscr{F}_{k,N}(x)=-\int_0^\infty\left(\frac{1}{1-e^{-t}}-\frac{1}{t}\right){}_N\textup{Li}_k\left(e^{-xt}\right)\, dt,
\end{align}
where ${}_N\textup{Li}_k(t)$ is defined in \eqref{lintN}, whereas Binet's formula \cite[p.~251]{ww} for Re$(x)>0$, namely
\begin{equation*}
\psi(x)+\frac{1}{2x}-\log x=-\int_{0}^{\infty}\frac{2t}{(t^2+x^2)(e^{2\pi t}-1)}\, dt,
\end{equation*}
yields a representation for $\mathscr{F}_{k, N}(x)$ in terms of a generalized Lambert series for Re$(x)>0$:
\begin{equation}
\mathscr{F}_{k, N}(x)=-\int_{0}^{\infty}\left(\sum_{n=1}^{\infty}\frac{n^{-k}}{\exp(2\pi n^{N}t)-1}\right)\frac{2t}{t^2+x^2}\, dt-\frac{1}{2x}\zeta(k+N).\nonumber
\end{equation}

\section{Proofs of functional equations satisfied by $\mathscr{F}_{k, N}(x)$}\label{proofsfe}

\begin{proof}[Theorem \textup{\ref{hhfe}}][]
Employing the well-known functional equation
\begin{align}\label{functionalpsi}
\psi(x+1)= \psi(x)+\frac{1}{x}
\end{align} 
in the definition of $\mathscr{F}_{k,N}(x)$ in \eqref{fknx}, we get
\begin{align}\label{befkloo}
\mathscr{F}_{k,N}(x)
=\sum_{n=1}^\infty\frac{\psi(n^Nx+1)-\log(n^Nx)}{n^k} -\frac{1}{x}\zeta(k+N).
\end{align}
We need the following formula due to Kloosterman \cite[Equation 2.9.2]{titch}, which is valid in $0<c=\textup{Re}(z)<1$:
\begin{align}\label{kloo}
\frac{1}{2\pi i}\int_{(c)}\frac{-\pi\zeta(1-z)}{\sin(\pi z)}x^{-z}\, dz=\psi(x+1)-\log x,
\end{align}
where here and in the sequel, $\int_{(c)}$ denotes the line integral $\int_{c-i\infty}^{c+i\infty}$. 

Thus, invoking \eqref{kloo} in \eqref{befkloo}, we get for $0<c=\textup{Re}(z)<1$,
\begin{align*}
\mathscr{F}_{k,N}(x)&=\frac{1}{2\pi i}\sum_{n=1}^\infty\frac{1}{n^k}\int_{(c)}\frac{-\pi\zeta(1-z)}{\sin (\pi z)} \left(n^{N}x \right)^{-z}\, dz-\frac{1}{x}\zeta(k+N) \nonumber\\
&=\frac{1}{2\pi i}\int_{(c)}\frac{-\pi\zeta(1-z)\zeta\left(k+Nz \right)}{\sin (\pi z)} x^{-z}\, dz-\frac{1}{x}\zeta(k+N),
\end{align*}
by interchanging the order of summation and integration. Employing the change of variable $z=-\frac{s}{N}- \frac{k-1}{N}$ in the above equation, we see that for $-N-k+1<c'=\textup{Re}(s)<1-k,$
\begin{align}\label{lineintegralkN}
\mathscr{F}_{k,N}(x)&=\frac{1}{N}x^{\frac{k-1}{N}}I-\frac{1}{x}\zeta(k+N), 
\end{align}
where, 
\begin{align}
I=I(k,N,x):=\frac{1}{2\pi i}\int_{(c')}\frac{\pi\zeta\left(1+\frac{s+k-1}{N}\right)\zeta(1-s)}{\sin \left(\frac{\pi}{N}(s+k-1) \right)} x^{\frac{s}{N}}\, ds.\nonumber
\end{align}
Now construct a rectangular contour formed by the line segments $[c'-iT, c'+iT]$, $[c'+iT, c''+iT]$, $[c''+iT, c''-iT]$, $[c''-iT, c'-iT]$ where $0<c''<N-k+1$. Since we wish to employ \eqref{kloo} again after shifting the line of integration from Re$(s)=c'$ to Re$(s)=c''$, we also need $0<c''<1$. In any case, we must have $N-k+1>0$, or, in other words, $k\leq N$. Note that in the process of shifting the line of integration, we encounter a double pole of the integrand at $s=1-k$ and a simple pole at $s=0$. 

Here, and throughout the paper, $R_a$ denotes the residue of the associated integrand at its pole at $a$. Stirling's formula in the vertical strip $p\leq\sigma\leq q$ reads \cite[p.~224]{cop}
\begin{equation}\label{strivert}
  |\Gamma(s)|=\sqrt{2\pi}|t|^{\sigma-\frac{1}{2}}e^{-\frac{1}{2}\pi |t|}\left(1+O\left(\frac{1}{|t|}\right)\right)
\end{equation}
as $|t|\to \infty$.

Then letting $T \to \i$, noting that the integrals over the horizontal segments thereby go to zero (due to \eqref{strivert}) and employing Cauchy's residue theorem, we have
\begin{align}\label{inter}
I=I_1-(R_0+ R_{1-k}), 
\end{align}
where
\begin{equation}
I_1=I_1(k,N,x):=\frac{1}{2\pi i}\int_{(c'')}\frac{\pi\zeta\left(1+\frac{s+k-1}{N}\right)\zeta(1-s)}{\sin \left(\frac{\pi}{N}(s+k-1) \right)} x^{\frac{s}{N}}\, ds,\nonumber
\end{equation}
and the residues $R_0$ and $R_{1-k}$ are easily seen to be
\begin{align}
R_0&= -\pi \frac{\zeta\left(1+\frac{k-1}{N}\right)}{\sin \left(\frac{\pi}{N}(k-1) \right)},\label{residues}\\
R_{1-k}&= -Nx^{\frac{1-k}{N}}\left(-\left(\gamma+\log x \right) \zeta(k)+N\zeta'(k)\right)\label{residues1}.
\end{align}
To evaluate $I_1$, we use the following well known result \cite[Lemma 4.3]{dgkm}
\begin{align}\label{Chebi}
\frac{1}{\sin(z)} =\frac{1}{\sin (N z)}\sum_{j=-(N-1)}^{(N-1)}{\vphantom{\sum}}''e^{ijz}
\end{align}
with $z= \pi (s+k-1)/N$ so as to obtain
\begin{align}\label{aftkloo}
I_1&= \sum_{j=-(N-1)}^{(N-1)}{\vphantom{\sum}}''e^{i\pi j (k-1)/N}\frac{1}{2\pi i}\int_{(c'')}\frac{\pi\zeta(1-s)\zeta\left(1+\frac{s+k-1}{N}\right)}{\sin \left(\pi(s+k-1) \right)}\left(\frac{e^{-\frac{i\pi j}{N}}}{x^{1/N}}\right)^{-s}\, ds\nonumber\\
&=(-1)^{k}\sum_{j=-(N-1)}^{(N-1)}{\vphantom{\sum}}''e^{i\pi j (k-1)/N}\sum_{n=1}^{\i}\frac{1}{n^{1+\frac{k-1}{N}}}\frac{1}{2\pi i}\int_{(c'')}\frac{-\pi\zeta(1-s)}{\sin \left(\pi s\right)} \left(\left(\frac{n}{x}\right)^{\frac{1}{N}}e^{-\frac{i\pi j}{N}}\right)^{-s}ds,
\end{align}
where in the last step we used the series representation of $\zeta\left(1+\frac{s+k-1}{N}\right)$ since $k\geq 1$ and $c''>0$ and interchanged the order of summation and integration which can be easily justified by standard arguments. Now invoke \eqref{kloo} again in \eqref{aftkloo} to obtain
\begin{align}\label{I1}
I_1&=(-1)^{k}\sum_{j=-(N-1)}^{(N-1)}{\vphantom{\sum}}''e^{i\pi j (k-1)/N}\sum_{n=1}^{\i}\frac{1}{n^{1+\frac{k-1}{N}}}\left\{\psi\left(\left(\frac{n}{x}\right)^{\frac{1}{N}}e^{-\frac{i\pi j}{N}}+1 \right)-\log \left(\left(\frac{n}{x}\right)^{\frac{1}{N}}e^{-\frac{i\pi j}{N}} \right) \right\} \nonumber\\
&=(-1)^{k}\sum_{j=-(N-1)}^{(N-1)}{\vphantom{\sum}}''e^{i\pi j (k-1)/N}\sum_{n=1}^{\i}\frac{1}{n^{1+\frac{k-1}{N}}}\left\{\psi\left(\left(\frac{n}{x}\right)^{\frac{1}{N}}e^{-\frac{i\pi j}{N}} \right)-\log \left(\left(\frac{n}{x}\right)^{\frac{1}{N}}e^{-\frac{i\pi j}{N}} \right) \right\}\nonumber\\
&\quad+(-1)^kx^{1/N}\zeta\left( 1+\frac{k}{N}\right)\sum_{j=-(N-1)}^{(N-1)}{\vphantom{\sum}}''e^{i\pi j k/N},
\end{align}
where we used \eqref{functionalpsi} in the last step. From \eqref{Chebi} and the fact that $1<k\leq N$, we have
\begin{align}\label{chebi1}
\sum_{j=-(N-1)}^{(N-1)}{\vphantom{\sum}}''e^{i\pi j k/N}=\left\{
	\begin{array}{ll}
		(-1)^{N+1}N & \mbox{if } k=N \\
		0 & \mbox{if } 1<k<N.
	\end{array}
\right.
\end{align}
Using \eqref{chebi1} in \eqref{I1} and recalling the definition of \eqref{fknx}, we see that
\begin{align}\label{I1a}
I_1=(-1)^{k}\sum_{j=-(N-1)}^{(N-1)}{\vphantom{\sum}}''e^{i\pi j (k-1)/N}\mathscr{F}_{\frac{N+k-1}{N},\frac{1}{N}}\left(\frac{e^{-\frac{i\pi j}{N}}}{x^{1/N}}\right)+N\mathscr{B}(k, N, x),
\end{align}
where $\mathscr{B}(k, N, x)$ is defined in \eqref{bknx}. Substituting \eqref{residues}, \eqref{residues1} and \eqref{I1a} in \eqref{inter}, we see that
\begin{align}\label{inter1}
I&=(-1)^{k}\sum_{j=-(N-1)}^{(N-1)}{\vphantom{\sum}}''e^{i\pi j (k-1)/N}\mathscr{F}_{\frac{N+k-1}{N},\frac{1}{N}}\left(\frac{e^{-\frac{i\pi j}{N}}}{x^{1/N}}\right)+N\mathscr{B}(k, N, x)+\pi \frac{\zeta\left(1+\frac{k-1}{N}\right)}{\sin \left(\frac{\pi}{N}(k-1) \right)}\nonumber\\
&\quad+Nx^{\frac{1-k}{N}}\left(-\left(\gamma+\log x \right) \zeta(k)+N\zeta'(k)\right).
\end{align}
The proof now follows from substituting \eqref{inter1} in \eqref{lineintegralkN} and multiplying both sides of the resulting equation by $x^{\frac{1-k}{N}}$.
\end{proof}

\begin{proof}[Theorem \textup{\ref{genherglotzNodd}}][]
The fact that $k$ and $N$ are odd is used several times in the proof without mention.

\textbf{Case 1}: $\frac{1-k}{N}\neq-2\left\lfloor\frac{k}{2N}\right\rfloor$.\\
Let $J$ denote the left-hand side of \eqref{thmnotequaleqn}. Using \eqref{fknx}, we rewrite $J$ in the form
\begin{align}\label{zq1}
J&=-2\sum_{j=-\left(N-1\right)}^{N-1}{\vphantom{\sum}}''i^j\exp\left(-\frac{i\pi(k+N-1)j}{2N}\right)\sum_{n=1}^\infty\frac{1}{n^{\frac{k+N-1}{N}}}\Bigg\{\log\left(\left(\frac{2\pi n}{x}\right)^{\frac{1}{N}}e^{\frac{i\pi j}{2N}}\right)\nonumber\\
&\qquad-\frac{1}{2}\left[\psi\left(i\left(\frac{2\pi n}{x}\right)^{\frac{1}{N}}e^{\frac{i\pi j}{2N}}\right)+\psi\left(-i\left(\frac{2\pi n}{x}\right)^{\frac{1}{N}}e^{\frac{i\pi j}{2N}}\right)\right]\Bigg\}.
\end{align}
Define 
\begin{align}
X_{\ell,n}:=&2\pi \ell\left(\frac{2\pi n}{x}\right)^{\frac{1}{N}}\label{xmn}\\
X_{\ell,n,j}^*:=&2\pi \ell\left(\frac{2\pi n}{x}\right)^{\frac{1}{N}}e^{\frac{i\pi j}{2N}}=X_{\ell,n}e^{\frac{i\pi j}{2N}}\label{xmn*}
\end{align}
Employing \eqref{raabesumeqn} in \eqref{zq1}, we observe that
\begin{align}\label{zq2}
J&=-4\sum_{j=-\left(N-1\right)}^{N-1}{\vphantom{\sum}}''i^j\exp\left(-\frac{i\pi(k+N-1)j}{2N}\right)\sum_{n=1}^\infty\frac{1}{n^{\frac{k+N-1}{N}}}\sum_{\ell=1}^\infty\int_0^\infty\frac{t\cos t}{t^2+{X_{\ell,n,j}^{*}}^2}\ dt\nonumber\\
&=-4\sum_{j=-\left(N-1\right)}^{N-1}{\vphantom{\sum}}''i^j\exp\left(-\frac{i\pi(k+N-1)j}{2N}\right)\frac{\pi}{2}\sum_{n, \ell=1}^\infty\frac{1}{n^{\frac{k+N-1}{N}}}\frac{1}{2\pi i}\int_{(c_1)}\frac{\Gamma(s_1)}{\tan\left(\frac{\pi s_1}{2}\right)}\left(X_{\ell,n,j}^*\right)^{-s_1}\ ds_1,
\end{align}
where $1<c_1:=\mathrm{Re}(s_1)<2$. The last step follows from the following result \cite[Lemma 4.1]{dgkm}
\begin{align}\label{nagoyalem}
\frac{1}{2\pi i}\int_{(c_1)}\frac{\Gamma(s_1)}{\tan(\frac{\pi s_1}{2})}u^{-s_1}\, ds_1&=\frac{2}{\pi}\int_0^\infty \frac{t\cos t}{u^2+t^2}\ dt,
\end{align}
which is valid for $ 0< \mathrm{Re}(s_1)=c_1 < 2 $ and $ \mathrm{Re}(u) >0$ \footnote{We take $1<c_1=\textup{Re}(s_1)<2$ in the integral in \eqref{zq2} and not $0<c_1=\textup{Re}(s_1)<2$ so that we are able to use the series representation of $\zeta(s)$ while evaluating the integral in \eqref{beforetwicefunct}}.

From \eqref{xmn}, \eqref{xmn*} and \eqref{zq2} and interchanging the order of summation and integration, we have
\begin{align}\label{zq3}
J=\frac{-1}{i}\sum_{n, \ell=1}^\infty\frac{1}{n^{\frac{k+N-1}{N}}}\int_{(c_1)}\frac{\Gamma(s_1)}{\tan\left(\frac{\pi s_1}{2}\right)}\sum_{j=-\left(N-1\right)}^{N-1}{\vphantom{\sum}}''i^j\exp\left(-\frac{i\pi j}{2N}(s_1+k+N-1)\right)X_{\ell,n}^{-s_1}\ ds_1.
\end{align}
From \cite[Equation~(4.2)]{dgkm}, for $N$ odd,
\begin{align}\label{dgkmdesh}
\frac{\cos(Nz)}{\cos(z)}=(-1)^{\frac{N-1}{2}}\sum_{j=-(N-1)}^{N-1}{\vphantom{\sum}}''i^j \exp({-ijz}).
\end{align}
Using \eqref{dgkmdesh} with $z=\pi(s_1+k+N-1)/(2N)$ in \eqref{zq3}, we observe that
\begin{align}
J&=\frac{-1}{i}\sum_{n, \ell=1}^\infty\frac{1}{n^{\frac{k+N-1}{N}}}\int_{(c_1)}\frac{\Gamma(s_1)}{\tan\left(\frac{\pi s_1}{2}\right)}(-1)^{\frac{N-1}{2}}\frac{\cos\left(\frac{\pi}{2}(s_1+k+N-1)\right)}{\cos\left(\frac{\pi}{2}\left(\frac{s_1+k+N-1}{N}\right)\right)}X_{\ell,n}^{-s_1}\ ds_1\nonumber\\
&=\frac{(-1)^{\frac{k-1}{2}}}{2i}\sum_{n, \ell=1}^\infty\frac{1}{n^{\frac{k+N-1}{N}}}\int_{(c_1)}\frac{\Gamma(s_1)\cos\left(\frac{\pi s_1}{2}\right)}{\cos\left(\frac{\pi}{2}\left(\frac{s_1+k+N-1}{N}\right)\right)}X_{\ell,n}^{-s_1}\ ds_1.\nonumber
\end{align}
Performing the change of variable $s_1=Ns-k-N+1$ in the above equation, we see that for $\frac{k+N}{N}<c<\frac{k+N+1}{N}$,
\begin{align}\label{beforetwicefunct}
J&=-\frac{(-1)^{\frac{N+1}{2}}N}{i}\sum_{n, \ell=1}^\infty\frac{1}{n^{\frac{k+N-1}{N}}}\int_{(c_1)}\frac{\Gamma(Ns-k-N+1)\sin\left(\frac{\pi Ns}{2}\right)}{\cos\left(\frac{\pi s}{2}\right)}X_{\ell,n}^{-Ns+k+N-1}\ ds\nonumber\\
&=\frac{-N}{i}(-1)^{\frac{N+1}{2}}\left(2\pi\left(\frac{2\pi}{x}\right)^{\frac{1}{N}}\right)^{k+N-1}\nonumber\\
&\quad\times\int_{(c)}\frac{\Gamma(Ns-k-N+1)\sin\left(\frac{\pi Ns}{2}\right)}{\cos\left(\frac{\pi s}{2}\right)}\zeta(s)\zeta(Ns-k-N+1)\left(2\pi\left(\frac{2\pi}{x}\right)^{\frac{1}{N}}\right)^{-Ns}\, ds,
\end{align}
where in the last step, we used the definition of $X_{\ell, n}$ in \eqref{xmn} and the series representations for $\zeta(s)$ and $\zeta(Ns-k-N+1)$.

Now use the asymmetric form of the functional equation for $\zeta(s)$ \cite[p.~13, Equation (2.1.1)]{titch}, namely,
\begin{equation}\label{zetaasym}
\zeta(s)=2^{s}\pi^{s-1}\G(1-s)\zeta(1-s)\sin\left(\frac{\pi s}{2}\right),
\end{equation}
 in \eqref{beforetwicefunct} to deduce that
\begin{align}\label{applicablelater}
J&=\frac{-N}{2\pi i}(-1)^{\frac{k+1}{2}}\left(\frac{2\pi}{x}\right)^{\frac{k+N-1}{N}}\int_{(c)}\zeta(k+N-Ns)\tan\left(\frac{\pi s}{2}\right)\Gamma(1-s)\zeta(1-s)x^s\ ds.
\end{align}
If we replace $s$ by $1-s$ in the above equation then for $-\frac{k+1}{N}<c_2<-\frac{k}{N}$, we get
\begin{align}\label{bcauchy}
J&=\frac{-N}{2\pi i}(-1)^{\frac{k+1}{2}}\left(\frac{2\pi}{x}\right)^{\frac{k+N-1}{N}}\int_{(c_2)}\zeta(k+Ns)\frac{\cos\left(\frac{\pi s}{2}\right)}{\sin\left(\frac{\pi s}{2}\right)}\Gamma(s)\zeta(s)x^{1-s}\ ds.
\end{align}
We wish to use the series representation for $\zeta(s)$ in \eqref{bcauchy} for which we transform it by shifting the line of integration to $1<d=\textup{Re}(s)<2$. Consider the contour formed by the lines $[c_2-iT,c_2+iT],\ [c_2+iT,d+iT],\ [d+iT,d-iT]$ and $[d-iT,c_2-iT]$. Note that inside this contour the integrand has\\
(1) a simple pole at $s=\frac{1-k}{N}$ (due to $\zeta(k+Ns)$);\\
(2) a double order pole at $s=0$ (due to $\sin\left(\frac{\pi s}{2}\right)$ and $\Gamma(s)$);\\
(3) simple poles at $s=-2j$ due to the zeros of $\sin\left(\frac{\pi s}{2}\right),\ 1\leq j\leq \left\lfloor\frac{k}{2N}\right\rfloor$ (because $-2j\geq-\frac{k}{N}$ which implies that $j\leq \left\lfloor\frac{k}{2N}\right\rfloor$). This is because the poles of $\Gamma(s)$ at $s=-j,\ 1\leq j\leq \frac{k}{N}$ are canceled by the zeros of $\zeta(s)$ and  $\cos\left(\frac{\pi s}{2}\right)$ at negative even and odd integers respectively.\\
Note that using \eqref{strivert}, one can see that the integrals over the horizontal segments of the contour tend to zero as the height $T\to\infty$. Therefore invoking Cauchy's theorem, we deduce that\\
\begin{align}\label{cauchyN}
&\frac{1}{2\pi i}\int_{(c_2)}\zeta(k+Ns)\cot\left(\frac{\pi s}{2}\right)\Gamma(s)\zeta(s)x^{1-s}\ ds\nonumber\\
&=-R_0-R_{\frac{1-k}{N}}-\sum_{j=1}^{\left\lfloor\frac{k}{2N}\right\rfloor}R_{-2j}+\frac{1}{2\pi i}\int_{(d)}\zeta(k+Ns)\cot\left(\frac{\pi s}{2}\right)\Gamma(s)\zeta(s)x^{1-s}\ ds.
\end{align}
The residues in \eqref{cauchyN} are calculated to be
\begin{align}\label{zeroN}
R_{0}&=\frac{x}{\pi}\left\{\left(\gamma-\log\left(\frac{2\pi}{x}\right)\right)\zeta(k)-N\zeta'(k)\right\},\nonumber\\
R_{-2j}&=2(-1)^j(2\pi)^{-2j-1}\zeta(k-2Nj)\zeta(1+2j)x^{1+2j},\nonumber\\
R_{\frac{1-k}{N}}&=\frac{\pi}{N}\left(\frac{x}{2\pi}\right)^{\frac{k+N-1}{N}}\frac{\zeta\left(1+\frac{k-1}{N}\right)}{\sin\left(\frac{\pi}{2}\left(\frac{1-k}{N}\right)\right)}.
\end{align}
Also, the integral on the right-hand side of \eqref{cauchyN} is evaluated to
\begin{align}\label{intN}
&\frac{1}{2\pi i}\int_{(d)}\zeta(k+Ns)\cot\left(\frac{\pi s}{2}\right)\Gamma(s)\zeta(s)x^{1-s}\ ds\nonumber\\
&=x\sum_{n, \ell=1}^\infty\frac{1}{\ell^{k}}\frac{1}{2\pi i}\int_{(d)}\frac{\Gamma(s)}{\tan\left(\frac{\pi s}{2}\right)}(x\ell^Nn)^{-s}\ ds\nonumber\\
&=\frac{2x}{\pi}\sum_{\ell=1}^\infty\frac{1}{\ell^{k}}\sum_{n=1}^\infty\int_0^\infty\frac{t\cos t}{(t^2+(x\ell^{N}n)^2)}\ dt\nonumber\\
&=\frac{x}{\pi}\sum_{\ell=1}^\infty\frac{1}{\ell^{k}}\left\{\log\left(\frac{\ell^Nx}{2\pi}\right)-\frac{1}{2}\left(\psi\left(\frac{i\ell^Nx}{2\pi}\right)+\left(\frac{-i\ell^Nx}{2\pi}\right)\right)\right\},
\end{align}
where in the second step we used \eqref{nagoyalem} and in the last step used \eqref{raabesumeqn}.

Now substitute \eqref{zeroN} and \eqref{intN} in \eqref{cauchyN} so as to derive
{\allowdisplaybreaks\begin{align}\label{evline}
&\frac{1}{2\pi i}\int_{(c)}\zeta(k+Ns)\cot\left(\frac{\pi s}{2}\right)\Gamma(s)\zeta(s)x^{1-s}\ ds\nonumber\\
&=\frac{x}{\pi}\sum_{\ell=1}^\infty\frac{1}{\ell^{k}}\left\{\log\left(\frac{\ell^Nx}{2\pi}\right)-\frac{1}{2}\left(\psi\left(\frac{i\ell^Nx}{2\pi}\right)+\left(\frac{-i\ell^Nx}{2\pi}\right)\right)\right\}\nonumber\\
&\quad-\frac{x}{\pi}\left\{\left(\gamma-\log\left(\frac{2\pi}{x}\right)\right)\zeta(k)-N\zeta'(k)\right\}\nonumber\\
&\quad-2\sum_{j=1}^{\left\lfloor\frac{k}{2N}\right\rfloor}(-1)^j(2\pi)^{-2j-1}\zeta(k-2Nj)\zeta(1+2j)x^{1+2j}\nonumber\\
&\quad-\frac{\pi}{N}\left(\frac{x}{2\pi}\right)^{\frac{k+N-1}{N}}\sec\left(\frac{\pi}{2}\left(\frac{k+N-1}{N}\right)\right)\zeta\left(\frac{k+N-1}{N}\right).
\end{align}}
Finally substitute \eqref{evline} in \eqref{bcauchy} and use the definitions of $\mathscr{F}_{k,N}(x)$ and $\mathscr{C}(k, N, x)$ from \eqref{fknx} and \eqref{cknx1} respectively to arrive at \eqref{thmnotequaleqn}.

\textbf{Case 2}: $\frac{1-k}{N}=-2\left\lfloor\frac{k}{2N}\right\rfloor\neq0$.\\
The only change in this case is the contribution of the double order pole of the integrand in \eqref{bcauchy} at $s=(1-k)/N$ due to $\zeta(k+Ns)$ and $\sin\left(\frac{\pi s}{2}\right)$. Due to this, the residue at $(1-k)/N$ of the integrand of the integral on the left-hand side of \eqref{cauchyN} now becomes
\begin{equation}
\frac{2(-1)^{\frac{k-1}{2N}}}{N}\left(\frac{x}{2\pi}\right)^{\frac{N+k-1}{N}}\left[\left(\gamma+\log\left(\frac{2\pi}{x}\right)\right)\zeta\left(1+\frac{k-1}{N}\right)-\zeta'\left(1+\frac{k-1}{N}\right)\right].\nonumber
\end{equation}
\end{proof}

\subsection{Proof of the equivalence of Theorems \ref{hhfe} and \ref{genherglotzNodd} for odd natural numbers $k$ and $N$ such that $1<k\leq N$}\label{equivalence}

Assume $k$ and $N$ to be odd positive integers such that  $1<k\leq N$. Employing \eqref{hhfeeqn} once, with $x$ replaced by $\frac{ix}{2\pi}$ and then again with $x$ replaced by $-\frac{ix}{2\pi}$, and  then adding the resulting two equations, we arrive at
\begin{align}\label{lhse1234}
\mathscr{F}_{k,N}\left(\frac{ix}{2\pi}\right)+\mathscr{F}_{k,N}\left(-\frac{ix}{2\pi}\right)&=\left(\frac{x}{2\pi}\right)^{\frac{k-1}{N}}\left\{E_1+E_2+E_3\right\}+E_4,
\end{align}
where 
\begin{align}
E_1:&=\frac{2\pi}{N}\frac{\zeta\left(\frac{N+k-1}{N}\right)}{\sin\left(\frac{\pi(k-1)}{N}\right)}\cos\left(\frac{\pi}{2N}(k-1)\right)\label{E1}\\
E_2:&=i^{\frac{k-1}{N}}\mathscr{B}\left(k,N,\frac{ix}{2\pi}\right)+(-i)^{\frac{k-1}{N}}\mathscr{B}\left(k,N,-\frac{ix}{2\pi}\right)\nonumber\\
E_3:&=\frac{(-1)^k}{N}{\sum_{j=-(N-1)}^{N-1}}{\vphantom{\sum}}''e^{\frac{i\pi j(k-1)}{N}}\left[i^{\frac{k-1}{N}}\mathscr{F}_{\frac{k+N-1}{N},\frac{1}{N}}\left(\left(\frac{ix}{2\pi}\right)^{-\frac{1}{N}}e^{-\frac{i\pi j}{N}}\right)\right.\nonumber\\
&\quad\left.+(-i)^{\frac{k-1}{N}}\mathscr{F}_{\frac{k+N-1}{N},\frac{1}{N}}\left(\left(-\frac{ix}{2\pi}\right)^{-\frac{1}{N}}e^{-\frac{i\pi j}{N}}\right)\right]\label{E3}\\
E_4:&=-2\gamma\zeta(k)+2N\zeta'(k)-\zeta(k)\left(\log\left(\frac{ix}{2\pi}\right)+\log\left(-\frac{ix}{2\pi}\right)\right)\label{E4}.
\end{align}
Using the definition of $\mathscr{B}(k,N,x)$ from \eqref{bknx}, it is easy to see that $E_2=0$ for $1<k<N$, and that for $k=N$, 
\begin{align*}
E_2&=2(-1)^{k+N+1}\left(\frac{x}{2\pi}\right)^{\frac{1}{N}}\zeta\left(1+\frac{k}{N}\right)\cos\left(\frac{\pi k}{2N}\right).
\end{align*}
Since the cosine function vanishes for $k=N$, we conclude that $E_2=0$ in this case as well. Hence 
\begin{equation}\label{E2simpli}
E_2=0\hspace{2mm}\text{for}\hspace{2mm} 1<k\leq N. 
\end{equation}
We next simplify $E_3$. To that end, use the definition of $\mathscr{F}_{k,N}(x)$ from \eqref{fknx} in \eqref{E3} so that 
{\allowdisplaybreaks\begin{align}
E_3&=\frac{(-1)^k}{N}{\sum_{j=-(N-1)}^{N-1}}{\vphantom{\sum}}''e^{\frac{i\pi j(k-1)}{N}}\left[i^{\frac{k-1}{N}}\sum_{n=1}^\infty \frac{\psi\left(\left(\frac{ix}{2\pi n}\right)^{-\frac{1}{N}}e^{-\frac{i\pi j}{N}}\right)-\log\left(\left(\frac{ix}{2\pi n}\right)^{-\frac{1}{N}}e^{-\frac{i\pi j}{N}}\right)}{n^{\frac{k+N-1}{N}}}\right.\nonumber\\
&\quad\left.+(-i)^{\frac{k-1}{N}}\sum_{n=1}^\infty \frac{\psi\left(\left(-\frac{ix}{2\pi n}\right)^{-\frac{1}{N}}e^{-\frac{i\pi j}{N}}\right)-\log\left(\left(-\frac{ix}{2\pi n}\right)^{-\frac{1}{N}}e^{-\frac{i\pi j}{N}}\right)}{n^{\frac{k+N-1}{N}}}\right].\nonumber
\end{align}}
Now write both the digamma functions in the above equation in the form $\psi(x)=\psi(x+1)-1/x$ by using the functional equation \eqref{functionalpsi} and then use Kloosterman's formula \eqref{kloo} twice in the above equation so as to get, for $0<c=\textup{Re}(s)<1$,
\begin{align}\label{eq3}
E_3&=\frac{(-1)^k}{N}{\sum_{j=-(N-1)}^{N-1}}{\vphantom{\sum}}''e^{\frac{i\pi j(k-1)}{N}}\left[-\left(\frac{x}{2\pi}\right)^{\frac{1}{N}}e^{\frac{i\pi j}{N}}\left(e^{\frac{i\pi k}{2N}}+e^{-\frac{i\pi k}{2N}}\right)\zeta\left(1+\frac{k}{N}\right)\right.\nonumber\\
&\left.\quad+\frac{1}{2\pi i}\int_{(c)}\frac{-\pi\zeta(1-s)}{\sin(\pi s)}\left(\left(\frac{x}{2\pi n}\right)^{-\frac{1}{N}}e^{-\frac{i\pi j}{N}}\right)^{-s}\left(i^{\frac{s+k-1}{N}}+(-i)^{\frac{s+k-1}{N}}\right)\ ds\right]\nonumber\\
&=\frac{2(-1)^{k+1}}{N}\left(\frac{x}{2\pi}\right)^{\frac{1}{N}}\cos\left(\frac{\pi k}{2N}\right)\zeta\left(1+\frac{k}{N}\right){\sum_{j=-(N-1)}^{N-1}}{\vphantom{\sum}}''e^{\frac{i\pi jk}{N}}+\frac{2(-1)^k}{2\pi iN}{\sum_{j=-(N-1)}^{N-1}}{\vphantom{\sum}}''e^{\frac{i\pi j(k-1)}{N}}\nonumber\\
&\quad\times\int_{(c)}\frac{-\pi\zeta(1-s)}{\sin(\pi s)}\zeta\left(\frac{s+k+N-1}{N}\right)\cos\left(\frac{\pi(s+k-1)}{2N}\right)\left(\left(\frac{x}{2\pi}\right)^{-\frac{1}{N}}e^{-\frac{i\pi j}{N}}\right)^{-s}\ ds.
\end{align}
We now use \eqref{Chebi} in \eqref{eq3} and employ the change of variable $s=s_1-(k+N-1)$ in the integral on the right-hand side of \eqref{eq3} to arrive at, for $k+N-1<\lambda=\textup{Re}(s_1)<k+N$,
\begin{align}
E_3&=\frac{(-1)^{k+1}}{N}\left(\frac{x}{\pi}\right)^{\frac{1}{N}}\zeta\left(1+\frac{k}{N}\right)\frac{\sin(\pi k)}{\sin\left(\frac{\pi k}{2N}\right)}+\frac{2(-1)^k}{2\pi iN}\left(\frac{x}{2\pi}\right)^{-\frac{(k+N-1)}{N}}{\sum_{j=-(N-1)}^{N-1}}{\vphantom{\sum}}''e^{-i\pi j}\nonumber\\
&\quad\times \int_{(\lambda)}\frac{\pi\zeta(k+N-s_1)}{(-1)^{k+N}\sin(\pi s_1)}\zeta\left(\frac{s_1}{N}\right)\sin\left(\frac{\pi s_1}{2N}\right)\left(\left(\frac{x}{2\pi}\right)^{-\frac{1}{N}}e^{-\frac{i\pi j}{N}}\right)^{-s}\ ds_1.\nonumber
\end{align}
Note that the first term on the right-hand side of the above equation vanishes as the sine function is zero since $k$ is an integer. Also as $k$ and $N$ are odd positive integers therefore $(-1)^{k+N}=1$ and $j$ runs over even integers hence we have $e^{-i\pi j}=1$.  Therefore $E_3$ reduces to 
\begin{align}\label{eq4}
E_3&=2\frac{(-1)^k}{2\pi iN}\left(\frac{x}{2\pi}\right)^{-\frac{(k+N-1)}{N}}{\sum_{j=-(N-1)}^{N-1}}{\vphantom{\sum}}''\int_{(\lambda)}\frac{\pi\zeta(k+N-s_1)}{\sin(\pi s_1)}\zeta\left(\frac{s_1}{N}\right)\sin\left(\frac{\pi s_1}{2N}\right)\nonumber\\
&\quad\times\left(\left(\frac{x}{2\pi}\right)^{-\frac{1}{N}}e^{-\frac{i\pi j}{N}}\right)^{-s_1}\ ds_1.
\end{align} 
Employing the change of variable $s_1=Ns$ in the integral in \eqref{eq4}, we have, for $\frac{k+N-1}{N}<\lambda^*=\textup{Re}(s)<\frac{k+N}{N}$,
\begin{align}\label{e3simp}
E_3&=\frac{2(-1)^k}{2\pi i}\left(\frac{x}{2\pi}\right)^{-\frac{(k+N-1)}{N}}{\sum_{j=-(N-1)}^{N-1}}{\vphantom{\sum}}''\int_{(\lambda^*)}\frac{\pi\zeta(k+N-Ns)}{\sin(\pi Ns)}\zeta\left(s\right)\sin\left(\frac{\pi s}{2}\right)\left(\frac{xe^{i\pi j}}{2\pi}\right)^{s}\ ds\nonumber\\
&=\frac{2(-1)^k}{2\pi i}\left(\frac{x}{2\pi}\right)^{-\frac{(k+N-1)}{N}}\int_{(\lambda^*)}\frac{\zeta(k+N-Ns)}{\sin(\pi Ns)}\sin^2\left(\frac{\pi s}{2}\right)\Gamma(1-s)\zeta\left(1-s\right)x^{s}{\sum_{j=-(N-1)}^{N-1}}{\vphantom{\sum}}''e^{i\pi js}ds,
\end{align}
where in the last step we used \eqref{zetaasym}. Now use \eqref{Chebi} in \eqref{e3simp} to obtain
\begin{align}\label{eq5}
E_3&=\frac{(-1)^k}{2\pi i}\left(\frac{x}{2\pi}\right)^{-\frac{(k+N-1)}{N}}\int_{(\lambda^*)}\zeta(k+N-Ns)\Gamma(1-s)\zeta\left(1-s\right)\tan\left(\frac{\pi s}{2}\right)x^{s}ds.
\end{align}
From \eqref{zq1} and \eqref{applicablelater}, for $\frac{k+N}{N}<c<\frac{k+N+1}{N}$, we have
\begin{align}\label{differentev}
&\frac{1}{2\pi i}\int_{(c)}\zeta(k+N-Ns)\Gamma(1-s)\zeta\left(1-s\right)\tan\left(\frac{\pi s}{2}\right)x^{s}ds\nonumber\\
&=\frac{-(-1)^{\frac{k+1}{2}}}{N}\left(\frac{x}{2\pi}\right)^{\frac{k+N-1}{N}}{\sum_{j=-(N-1)}^{N-1}}{\vphantom{\sum}}''e^{-\frac{i\pi(k-1)j}{2N}}\left\{\mathscr{F}_{\frac{k+N-1}{N},\frac{1}{N}}\left(i\left(\frac{2\pi}{x}\right)^{\frac{1}{N}}e^{\frac{i\pi j}{2N}}\right)\right.\nonumber\\
&\left.\quad+\mathscr{F}_{\frac{k+N-1}{N},\frac{1}{N}}\left(-i\left(\frac{2\pi}{x}\right)^{\frac{1}{N}}e^{\frac{i\pi j}{2N}}\right)\right\}.
\end{align}
We want to employ \eqref{differentev} in \eqref{eq5}, for which, we need to shift the line of integration of \eqref{eq5} to $\frac{k+N}{N}<c<\frac{k+N+1}{N}$. Consider the contour formed by the line segments $[\lambda^*-iT,c-iT],\ [c-iT,c+iT],\ [c+iT,\lambda^*+iT]$ and $[\lambda^*+iT,\lambda^*-iT]$. Note that the integrand has no poles inside this contour. The condition $1<k\leq N$ implies $-1\leq1-\frac{k+N}{N}<\frac{-1}{N}$, ensures that the pole of the integrand of the left-hand side of \eqref{differentev} at $s=\frac{k+N}{N}$ does not lie inside the contour. Due to the same reason, $\zeta(1-s)$ also does not have a pole inside the contour. Therefore applying Cauchy's residue theorem to the integral in \eqref{eq5} and using the fact that the integral along the horizontal segments go to zero as the height $T$ of the contour tends to $\infty$, for $\frac{k+N}{N}<c<\frac{k+N+1}{N}$, we see that
\begin{align}\label{eq6}
E_3&=\frac{(-1)^k}{2\pi i}\left(\frac{x}{2\pi}\right)^{-\frac{(k+N-1)}{N}}\int_{(c)}\zeta(k+N-Ns)\Gamma(1-s)\zeta\left(1-s\right)\tan\left(\frac{\pi s}{2}\right)x^{s}ds\nonumber\\
&=\frac{(-1)^{\frac{k+1}{2}}}{N}{\sum_{j=-(N-1)}^{N-1}}{\vphantom{\sum}}''e^{-\frac{i\pi(k-1)j}{2N}}\left\{\mathscr{F}_{\frac{k+N-1}{N},\frac{1}{N}}\left(i\left(\frac{2\pi}{x}\right)^{\frac{1}{N}}e^{\frac{i\pi j}{2N}}\right)\right.\nonumber\\
&\left.\quad+\mathscr{F}_{\frac{k+N-1}{N},\frac{1}{N}}\left(-i\left(\frac{2\pi}{x}\right)^{\frac{1}{N}}e^{\frac{i\pi j}{2N}}\right)\right\},
\end{align}
where in the last step we used \eqref{differentev}.

Note that $k$ and $N$ odd and $1<k\leq N$ imply that $\frac{1-k}{N}\neq-2\left\lfloor\frac{k}{2N}\right\rfloor$, which, in turn, implies that the finite sum in \eqref{cknx1} is empty. Therefore, the equivalence of  Theorems \ref{hhfe} and \ref{genherglotzNodd} follows upon substituting \eqref{E1}, \eqref{E4}, \eqref{E2simpli} and \eqref{eq6} in \eqref{lhse1234} to arrive at \eqref{thmnotequaleqn}.

\begin{proof}[Theorem \textup{\ref{thmtriplepole}}][]
Let $k=1$ in \eqref{bcauchy}. The proof is almost similar to that of Theorem \ref{genherglotzNodd} except that the pole of the integrand in \eqref{bcauchy} at $0$ is of order three whence the corresponding residue becomes
\begin{align*}
R_0=\frac{-x}{2N\pi}\left\{-2\g^2N+\frac{\pi^{2}}{12}+2\g(N-1)\log\left(\frac{2\pi}{x}\right)+\log^{2}(2\pi)-\log\left(\frac{4\pi^2}{x}\right)\log(x)-2(N^2+1)\g_1\right\}.
\end{align*}
\end{proof}

\begin{proof}[Corollary \textup{\ref{trans2m+1}}][]
Let $N=1$ in Theorem \ref{genherglotzNodd} and use the definition of $\mathscr{F}_{k,N}(x)$ in \eqref{fknx}. Since $m\geq1$, we can separate the expressions involving logarithm. This results in a lot of simplification thereby resulting in \eqref{trans2m+1eqn}.
\end{proof}

\begin{proof}[Corollary \textup{\ref{modular}}][]
Let $m=1$ in Corollary \ref{trans2m+1} and notice that the finite sum on the right-hand side vanishes thereby giving \eqref{modulareqn}.
\end{proof}

\subsection{Proof of a relation between $\mathscr{F}_{k, N}(x)$ and a generalized polylogarithm}\label{relation}

\begin{proof}[Theorem \textup{\ref{zagintevalgen}}][]

Employing \eqref{intreprF} and repeatedly using the fact  $\displaystyle\sum_{n=1}^{\infty}a(2n)=\frac{1}{2}\sum_{n=1}^{\infty}(1+(-1)^n)a(n)$, we see that
\begin{align}
&\mathscr{F}_{k,N}(2^Nx)-\left(2^{k-1}+2^{1-k}\right)\mathscr{F}_{k,N}(x)+\mathscr{F}_{k,N}\left(\frac{x}{2^N}\right)\nonumber\\
&=-\int_0^\infty\left(\frac{1}{1-e^{-t}}-\frac{1}{t}\right)\left[\sum_{n=1}^\infty \frac{e^{-(2n)^Nxt}}{n^k}-\left(2^{k-1}+2^{1-k}\right)\sum_{n=1}^\infty\frac{e^{-n^Nxt}}{n^k}+\sum_{n=1}^\infty \frac{e^{-\left(n/2\right)^Nxt}}{n^k}\right]dt\nonumber\\
&=-\int_0^\infty\left(\frac{1}{1-e^{-t}}-\frac{1}{t}\right)\left[2^{k-1}\sum_{n=1}^\infty(1+(-1)^n)\frac{e^{-n^Nxt}}{n^k}-2^{k-1}\sum_{n=1}^\infty\frac{e^{-n^Nxt}}{n^k}\right.\nonumber\\
&\quad\left.-\sum_{n=1}^\infty\left(1+(-1)^n\right)\frac{e^{-(n/2)^Nxt}}{n^k}+\sum_{n=1}^\infty\frac{e^{-(n/2)^Nxt}}{n^k}\right]dt\nonumber\\
&=-\int_0^\infty\left(\frac{1}{1-e^{-t}}-\frac{1}{t}\right)\left[2^{k-1}\sum_{n=1}^\infty(-1)^n\frac{e^{-n^Nxt}}{n^k}-\sum_{n=1}^\infty(-1)^n\frac{e^{-(n/2)^Nxt}}{n^k}\right]\, dt\nonumber\\
&=-\int_0^\infty\left(\frac{1}{e^{t}-1}-\frac{1}{t}\right)\left[2^{k-1}\sum_{n=1}^\infty(-1)^n\frac{e^{-n^Nxt}}{n^k}-\sum_{n=1}^\infty(-1)^n\frac{e^{-(n/2)^Nxt}}{n^k}\right]\, dt\nonumber\\
&\quad-\int_0^\infty\left(2^{k-1}\sum_{n=1}^\infty(-1)^n\frac{e^{-n^Nxt}}{n^k}-\sum_{n=1}^\infty(-1)^n\frac{e^{-(n/2)^Nxt}}{n^k}\right)dt,\nonumber
\end{align}
where the last step was simplified using $\displaystyle\frac{1}{1-e^{-t}}=\frac{1}{e^{t}-1}+1$. Since $\displaystyle\int_{0}^{\infty}e^{-a t}\, dt=\frac{1}{a}$ for $a>0$, we have
\begin{align}
&\mathscr{F}_{k,N}(2^Nx)-\left(2^{k-1}+2^{1-k}\right)\mathscr{F}_{k,N}(x)+\mathscr{F}_{k,N}\left(\frac{x}{2^N}\right)\nonumber\\
&=-2^{k-1}\int_0^\infty\left(\frac{1}{e^{t}-1}-\frac{1}{t}\right)\sum_{n=1}^\infty(-1)^n\frac{e^{-n^Nxt}}{n^k}\, dt+2^{N}\int_0^\infty\left(\frac{1}{e^{2^{N}y}-1}-\frac{1}{2^{N}y}\right)\sum_{n=1}^\infty(-1)^n\frac{e^{-n^Nxy}}{n^k}\, dy\nonumber\\
&\quad-\frac{\left(2^{k-1}-2^N\right)}{x}\sum_{n=1}^\infty\frac{(-1)^n}{n^{k+N}},\nonumber
\end{align}
where, in the second integral, we employed the change of variable $t=2^{N}y$. Since $\sum_{n=1}^{\infty}(-1)^{n}n^{-s}=(2^{1-s}-1)\zeta(s)$ for Re$(s)>0$, we deduce that
\begin{align}\label{FJ3}
&\mathscr{F}_{k,N}(2^Nx)-(2^{k-1}+2^{1-k})\mathscr{F}_{k,N}(x)+\mathscr{F}_{k,N}\left(\frac{x}{2^N}\right)\nonumber\\
&=-\int_0^\infty\left(\frac{2^{k-1}}{e^{t}-1}-\frac{2^N}{e^{2^Nt}-1}-\frac{(2^{k-1}-1)}{t}\right)\sum_{n=1}^\infty(-1)^n\frac{e^{-n^Nxt}}{n^k}dt\nonumber\\
&\quad-(2^N+2^{-N}-2^{k-1}-2^{1-k})\frac{\zeta(k+N)}{x}.
\end{align}
Making change of variable $e^{-t}=u$ in the integral on the right-hand side of \eqref{FJ3}, using the definitions of ${}_N\textup{Li}_{s}(t)$ and $J_{k,N}(x)$ given in \eqref{lintN} and \eqref{defnJkN} respectively, using the elementary fact $(-1)^{n^{N}}=(-1)^{n}$ and simplifying, we arrive at \eqref{zagintevalgeneqn}. This completes the proof.

\end{proof}

The next result which gives a relation between an integral with polylogarithm in its integrand and the Vlasenko-Zagier higher Herglotz function, and is analogous to \eqref{zaginteval}, was missing in the literature.
\begin{proof}[Corollary \textup{\ref{N=1JK}}][]
Let $N=1$ in the Theorem \ref{zagintevalgen} and employ \eqref{fknxk>1}.
\end{proof}

\begin{proof}[Corollary \textup{\ref{N=1x1}}][]
Let $x=1$ in Corollary \ref{N=1JK}. This gives
\begin{align}\label{x=1cor2.9}
&\int_0^1\left(\frac{2^{k-1}}{u-1}-\frac{2u}{u^{2}-1}-\frac{2^{k-1}-1}{u\log u}\right)\mathrm{Li}_k(-u) \ du\nonumber\\
&=F_k(2)-\left(2^{k-1}+2^{1-k}\right)F_{k}(1)+F_{k}\left(\frac{1}{2}\right)+\left(2-2^{k-1}-2^{1-k}\right)\left(\zeta'(k)+\zeta(k+1)\right)+\frac{1}{2}\zeta(k+1).
\end{align}
Now let $k>1$ be an odd integer. We first calculate $F_k(1)$. Letting $x=1$ in \eqref{vz1f} and simplifying, we obtain
\begin{align}\label{fk1value}
F_k(1)=-\g\zeta(k)-\zeta(k+1)-\frac{1}{2}\sum_{r=2}^{k-1}(-1)^{r-1}\zeta(r)\zeta(k+1-r).
\end{align}
Also, employing \eqref{vz1f} again with $x=2$ and simplifying, we deduce that
\begin{align}\label{fkhalfvalue}
F_k\left(\frac{1}{2}\right)=-2^{1-k}F_k(2)-2(1+2^{-k-1})\zeta(k+1)-(1+2^{1-k})\g\zeta(k)-2^{1-k}\sum_{r=2}^{k-1}(-2)^{r-1}\zeta(r)\zeta(k+1-r).
\end{align}
Now, substitute \eqref{fk1value} and \eqref{fkhalfvalue} in \eqref{x=1cor2.9} so as to obtain \eqref{N=1x1eqn} upon simplification.
\end{proof}

\section{Proofs of transformations and asymptotic expansions of the extended higher Herglotz functions and generalized Lambert series}\label{proofsasym}
\begin{proof}[Theorem \textup{\ref{asymptotic}}][]
The formula \cite[p.~259, formula 6.3.18]{as}
\begin{align}
\psi(x) \sim \log x + \sum_{n=1}^{\i}\zeta(1-n)x^{-n},\nonumber
\end{align}
as $x\to \infty$, $|\arg$ $x|<\pi$ implies that
\begin{align}
\mathscr{F}_{k,N}(x) &\sim\sum_{n=1}^{\infty}\zeta(1-n)\zeta(k+nN)x^{-n}\label{required}\\
&=-\frac{1}{2x}\zeta(k+N)-\sum_{n=1}^{\infty}\frac{B_{2n}}{2n}\zeta(k+2nN)x^{-2n},\nonumber
\end{align}
where in the second step we used the facts that $\zeta(-2j)=0$ for $j\in\mathbb{N}$ and \cite[p.~266, Theorem 12.16]{apostol}
\begin{equation}\label{zeta-}
\zeta(-n)=\begin{cases}
-\frac{1}{2},\hspace{2mm}\text{if}\hspace{1mm}n=0,\\
-\frac{B_{n+1}}{n+1},\hspace{2mm}\text{if}\hspace{1mm}n\in\mathbb{N}.
\end{cases}
\end{equation}

We now establish \eqref{f1nx} using \eqref{hhfeeqn} and \eqref{required}. Upon using \eqref{required}, as $x \to 0$, we have 
{\allowdisplaybreaks\begin{align}\label{multisumj}
&\sum_{j=-(N-1)}^{(N-1)}{\vphantom{\sum}}''e^{i\pi j (k-1)/N}\mathscr{F}_{\frac{N+k-1}{N},\frac{1}{N}}\left(\left(\frac{e^{-i\pi j}}{x}\right)^{1/N}\right) \nonumber \\
& \sim \sum_{j=-(N-1)}^{(N-1)}{\vphantom{\sum}}''e^{i\pi j (k-1)/N}\sum_{n=1}^{\i}\frac{\zeta(1-n)\zeta\left(\frac{N+k+n-1}{N}\right)}{\left(\left(\frac{e^{-i\pi j}}{x}\right)^{1/N}\right)^n} \nonumber\\
&= \sum_{n=1}^{\i}\zeta(1-n)\zeta\left(\frac{N+k+n-1}{N}\right)x^{n/N}\sum_{j=-(N-1)}^{(N-1)}{\vphantom{\sum}}''e^{i\pi j (k+n-1)/N} \nonumber\\
&= \sum_{n=1}^{\i}\zeta(1-n)\zeta\left(\frac{N+k+n-1}{N}\right)x^{n/N}\sum_{j=-(N-1)}^{(N-1)}{\vphantom{\sum}}''e^{i\pi j (k+n-1)/N}.
\end{align}}
From \eqref{Chebi}, we know that if $m \in \mathbb{N}$,
\begin{align}\label{chebinew}
\sum_{j=-(N-1)}^{(N-1)}{\vphantom{\sum}}''e^{i\pi j(k+n-1)/N}=\left\{
	\begin{array}{ll}
		(-1)^{m(N+1)}N & \mbox{if } n=mN-k+1 \\
		0 & \mbox{if } n\neq mN-k+1 .
	\end{array}
\right.
\end{align}
Employing \eqref{chebinew} in \eqref{multisumj}, we see that, as $x \to 0$,
\begin{align}\label{5.5}
&\sum_{j=-(N-1)}^{(N-1)}{\vphantom{\sum}}''e^{i\pi j (k-1)/N}\mathscr{F}_{\frac{N+k-1}{N},\frac{1}{N}}\left(\left(\frac{e^{-i\pi j}}{x}\right)^{1/N}\right) \nonumber\\
&\sim  N\sum_{m\geq \frac{k}{N}}^{\i}(-1)^{m(N+1)}\zeta(k-mN)\zeta\left(1+m\right)x^{m+\frac{1-k}{N}}.
\end{align}
Now divide \eqref{hhfeeqn} by $x^{(1-k)/N}$, then take $x \to 0$ in the resulting identity to obtain \eqref{f1nx}. 

Equation \eqref{ff1nx} can be proved by letting $k=1$ in \eqref{5.5} and putting the resulting asymptotic expansion in \eqref{hhfeeqn1}.
\end{proof}

\begin{proof}[Theorem \textup{\ref{asymthm2.3}}][]
From \eqref{fknxinfty}, it is easy to observe that as $x \to \infty$,
\begin{align*}
\mathscr{F}_{k, N}\left(\frac{ix}{2\pi}\right)+\mathscr{F}_{k, N}\left(-\frac{ix}{2\pi}\right)
\sim\sum_{n=1}^{\infty}\frac{(-1)^{n+1}}{n}B_{2n}\zeta(k+2nN)\left(\frac{2\pi}{x}\right)^{2n}.
\end{align*}
To obtain \eqref{f+f} as $x\to 0$, we first rewrite  \eqref{thmnotequaleqn} as,
\begin{align}\label{000}
&\mathscr{F}_{k,N}\left(\tfrac{ix}{2\pi}\right)+\mathscr{F}_{k,N}\left(-\tfrac{ix}{2\pi}\right)=-\bigg\{2\left(\left(\gamma-\log\left(\tfrac{2\pi}{x}\right)\right)\zeta(k)-N\zeta'(k)\right)+\mathscr{C}(k, N, x)\bigg\}+\frac{(-1)^{\frac{k+1}{2}}}{N}\left(\frac{x}{2\pi}\right)^{\frac{k-1}{N}}\nonumber\\
&\times\sum_{j=-(N-1)}^{N-1}{\vphantom{\sum}}''\exp\left(\tfrac{i\pi j}{2}-\tfrac{i\pi(N+k-1)}{2N}\right)\left\{\mathscr{F}_{\frac{N+k-1}{N},\frac{1}{N}}\left(i\left(\tfrac{2\pi}{x}\right)^{\frac{1}{N}}e^{\frac{i\pi j}{2N}}\right)+\mathscr{F}_{\frac{N+k-1}{N},\frac{1}{N}}\left(-i\left(\tfrac{2\pi}{x}\right)^{\frac{1}{N}}e^{\frac{i\pi j}{2N}}\right)\right\}.
\end{align}
Using \eqref{required}, as $x \to 0$,
\begin{align}
&\sum_{j=-(N-1)}^{N-1}{\vphantom{\sum}}''\exp\left(\tfrac{i\pi j}{2}\right)\left\{\mathscr{F}_{\frac{N+k-1}{N},\frac{1}{N}}\left(i\left(\tfrac{2\pi}{x}\right)^{\frac{1}{N}}e^{\frac{i\pi j}{2N}}\right)+\mathscr{F}_{\frac{N+k-1}{N},\frac{1}{N}}\left(-i\left(\tfrac{2\pi}{x}\right)^{\frac{1}{N}}e^{\frac{i\pi j}{2N}}\right)\right\}\nonumber\\
&\qquad\sim2\sum_{j=-(N-1)}^{N-1}{\vphantom{\sum}}''\exp\left(\tfrac{i\pi j}{2}\right)\sum_{n=1}^{\infty}\zeta(1-n)\zeta\left(\frac{k+n+N-1}{N}\right)\left(\frac{x}{2\pi} \right)^{\frac{n}{N}}e^{\frac{-i\pi nj}{2N}}\cos\left(\frac{\pi n}{2} \right)\nonumber\\
&\qquad=2\sum_{n=1}^{\infty}(-1)^n\zeta(1-2n)\zeta\left(\frac{k+2n+N-1}{N}\right)\left(\frac{x}{2\pi} \right)^{\frac{2n}{N}}\sum_{j=-(N-1)}^{N-1}{\vphantom{\sum}}''\exp\left(\frac{i\pi j}{2}-\frac{i\pi nj}{N}\right),\nonumber
\end{align}
where in the last step the fact that cosine function vanishes for odd values of $n$. Using \eqref{Chebi} and \eqref{zeta-} in the above equation, we arrive at 
\begin{align}\label{12}
&\sum_{j=-(N-1)}^{N-1}{\vphantom{\sum}}''\exp\left(\tfrac{i\pi j}{2}\right)\left\{\mathscr{F}_{\frac{N+k-1}{N},\frac{1}{N}}\left(i\left(\tfrac{2\pi}{x}\right)^{\frac{1}{N}}e^{\frac{i\pi j}{2N}}\right)+\mathscr{F}_{\frac{N+k-1}{N},\frac{1}{N}}\left(-i\left(\tfrac{2\pi}{x}\right)^{\frac{1}{N}}e^{\frac{i\pi j}{2N}}\right)\right\}\nonumber\\
&\qquad=-\sum_{n=1}^{\infty}\frac{B_{2n}}{n}\zeta\left(\frac{k+2n+N-1}{N}\right)\left(\frac{x}{2\pi} \right)^{\frac{2n}{N}}\frac{\sin\left(\frac{\pi N}{2}\right)}{\cos\left(\frac{\pi n}{N}\right)}.
\end{align}
Hence employ \eqref{12} in \eqref{000} to obtain \eqref{f+f}.

To obtain \eqref{f+f1}, let $k=1$ in \eqref{12} and then substitute the resultant in Theorem \ref{thmtriplepole}.
\end{proof}

\begin{proof}[Theorem \textup{\ref{gencompanion}}][]
We need Theorem 2.4 from \cite{dgkm} given below.\\

\textit{Let $0<a\leq 1$, let $N$ be an odd positive integer and $\alpha,\beta>0$ such that $\alpha\beta^{N}=\pi^{N+1}$. Then for any positive integer $m$,
\begin{align}\label{zetageneqna}
&\alpha^{-\frac{2Nm}{N+1}}\bigg(\left(a-\frac{1}{2}\right)\zeta(2Nm+1)+\sum_{j=1}^{m-1}\frac{B_{2j+1}(a)}{(2j+1)!}\zeta(2Nm+1-2jN)(2^N\alpha)^{2j}\nonumber\\
&\qquad\qquad+\sum_{n=1}^{\infty}\frac{n^{-2Nm-1}\textup{exp}\left(-a(2n)^{N}\alpha\right)}{1-\textup{exp}\left(-(2n)^{N}\alpha\right)}\bigg)\nonumber\\
&=\left(-\beta^{\frac{2N}{N+1}}\right)^{-m}\frac{2^{2m(N-1)}}{N}\bigg[\frac{(-1)^{m+1}(2\pi)^{2m}B_{2m+1}(a)N \gamma}{(2m+1)!}+\frac{1}{2}\sum_{n=1}^{\infty}\frac{\cos(2\pi na)}{n^{2m+1}}\nonumber\\
&\quad+(-1)^{\frac{N+3}{2}}\sum_{j=\frac{-(N-1)}{2}}^{\frac{N-1}{2}}(-1)^{j}\bigg\{\sum_{n=1}^{\infty}\frac{n^{-2m-1}\cos(2\pi na)}{\textup{exp}\left((2n)^{\frac{1}{N}}\beta e^{\frac{i\pi j}{N}}\right)-1}\nonumber\\
&\quad+\frac{(-1)^{j+\frac{N+3}{2}}}{2\pi}\sum_{n=1}^{\infty}\frac{\sin(2\pi na)}{n^{2m+1}}\left(\psi\left(\tfrac{i\beta}{2\pi}(2n)^{\frac{1}{N}} e^{\frac{i\pi j}{N}}\right)+\psi\left(\tfrac{-i\beta}{2\pi}(2n)^{\frac{1}{N}} e^{\frac{i\pi j}{N}}\right)\right)\bigg\}\bigg]\nonumber\\
&\quad+(-1)^{m+\frac{N+3}{2}}2^{2Nm}\sum_{j=0}^{\left\lfloor\frac{N+1}{2N}+m\right\rfloor}\frac{(-1)^jB_{2j}(a)B_{N+1+2N(m-j)}}{(2j)!(N+1+2N(m-j))!}\alpha^{\frac{2j}{N+1}}\beta^{N+\frac{2N^2(m-j)}{N+1}}.
\end{align}}

The main idea is to differentiate both sides of the above identity with respect to $a$ and then let $a=1$ while using the facts $\frac{d}{da}B_{2j}(a)=2jB_{2j-1}(a)$, $B_{n}(1)=(-1)^nB_n$, $B_{2 m }=\frac{2(-1)^{m+1}(2 m)!\zeta(2 m )}{(2\pi)^{2m}}$ and $\alpha\beta^{N}=\pi^{N+1}$. Inducting the $\zeta(2Nm+1)$ from the resulting identity in the finite sum on the resulting left-hand side as its $j=0$ term, we see that
\begin{align*}
&\alpha^{-\frac{2Nm}{N+1}}\left[\sum_{j=0}^{m-1}\frac{B_{2j+1}(a)}{(2j+1)!}\zeta(2Nm+1-2jN)(2^N\alpha)^{2j}-2^{N}\a\sum_{n=1}^{\infty}\frac{n^{-2Nm-1+N}}{\exp((2n)^{N}\a)-1}\right]\nonumber\\
&=\left(-\beta^{\frac{2N}{N+1}}\right)^{-m}\frac{2^{2m(N-1)}}{N}\bigg[2N\g\zeta(2m)+\sum_{j=\frac{-(N-1)}{2}}^{\frac{N-1}{2}}\sum_{n=1}^{\infty}\frac{1}{n^{2m}}\left(\psi\left(\frac{i\beta}{2\pi}(2n)^{\frac{1}{N}} e^{\frac{i\pi j}{N}}\right)\right.\nonumber\\
&\quad\left.+\psi\left(\frac{-i\beta}{2\pi}(2n)^{\frac{1}{N}} e^{\frac{i\pi j}{N}}\right)\right)\bigg]+2^{N-1}\a^{1-\frac{2Nm}{N+1}}\zeta(2Nm+1-N).
\end{align*}
The result now follows by dividing both sides of the above equation by $2^N\sqrt{\a}$, using the fact $\alpha\beta^{N}=\pi^{N+1}$ and rearranging the resulting equation.
\end{proof}

\begin{proof}[Theorem \textup{\ref{gendivdkk}}][]
Let $0<a\leq1$ and $N$ be an odd positive integer. Let $\alpha,\ \beta>0$ be such that $\alpha\beta=\pi^{N+1}$. Then from \cite[Theorem 2.7]{dgkm}, we have
\begin{align}\label{limab}
&\sum_{n=1}^{\infty}\frac{\textup{exp}(-a(2n)^{N}\a)}{n(1-\textup{exp}(-(2n)^{N}\a))}-\frac{1}{N}(-1)^{\frac{N+3}{2}}\sum_{j=-\frac{(N-1)}{2}}^{\frac{N-1}{2}}(-1)^{j}\bigg(\sum_{n=1}^{\infty}\frac{\cos(2\pi na)}{n\left(\textup{exp}\left((2n)^{\frac{1}{N}}\b e^{\frac{i\pi j}{N}}\right)-1\right)}\nonumber\\
&+\frac{(-1)^{j+\frac{N+1}{2}}}{\pi}\sum_{n=1}^{\infty}\frac{\sin(2\pi na)}{n}\left\{\log\left(\tfrac{\beta}{2\pi}(2n)^{\frac{1}{N}} e^{\frac{i\pi j}{N}}\right)-\frac{1}{2}\left(\psi\left(\tfrac{i\beta}{2\pi}(2n)^{\frac{1}{N}} e^{\frac{i\pi j}{N}}\right)+\psi\left(\tfrac{-i\beta}{2\pi}(2n)^{\frac{1}{N}} e^{\frac{i\pi j}{N}}\right)\right)\right\}\bigg)\nonumber\\
&=\frac{1}{N}\left((a-1)\log(2\pi)+\log\G(a)\right)+\left(a-\frac{1}{2}\right)\left\{\frac{(N-1)(\log 2-\g)}{N}+\frac{\log(\a/\b)}{N+1}\right\}\nonumber\\
&\quad+(-1)^{\frac{N+3}{2}}\sum_{j=0}^{\left\lfloor\frac{N+1}{2N}\right\rfloor}\frac{(-1)^jB_{2j}(a)B_{N+1-2Nj}}{(2j)!(N+1-2Nj)!}\a^{\frac{2j}{N+1}}\b^{N-\frac{2N^2j}{N+1}}.
\end{align}
Differentiate both sides of \eqref{limab} with respect to $a$, let $a=1$ in the resulting identity, and then divide both sides by $-2^N\alpha$ to arrive at \eqref{gendivdkkeqnalphabeta} upon simplification. 

\end{proof}

\begin{proof}[Theorem \textup{\ref{lambert asymptotic}}][]
Using \eqref{fknx}, we rewrite \eqref{gencompanioneqn} as 
\begin{align}\label{lambert}
\sum_{n=1}^{\i}&\frac{n^{-2Nm-1+N}}{e^{(2n)^N\a}-1}= -\frac{1}{2}\zeta(2Nm+1-N)+\frac{(-1)^{m+1}2^{2m(N-1)}\b^{N/2}}{N\pi^{\frac{N+1}{2}}2^N\sqrt{\a}}\frac{}{}\left(\frac{\b}{\a}\right)^{-\frac{2Nm}{N+1}}\nonumber\\
&\times\sum_{j=-\frac{(N-1)}{2}}^{\frac{(N-1)}{2}}\left(\mathscr{F}_{2m,\frac{1}{N}}\left( \frac{i \b}{2\pi}2^{1/N}e^{\frac{i\pi j}{N}}\right)+\mathscr{F}_{2m,\frac{1}{N}}\left( -\frac{i \b}{2\pi}2^{1/N}e^{\frac{i\pi j}{N}}\right)\right)+\a^{\frac{2Nm}{N+1}-\frac{1}{2}}\mathscr{E}_{N}(m,\a,\b),
\end{align}
where $\mathscr{E}_{N}(m,\a,\b)$ is defined by
\begin{align}\label{Dm}
&\mathscr{E}_{N}(m,\a, \b):=\frac{2^{(2m-1)(N-1)}}{N\pi^{\frac{N+1}{2}}}(-1)^{m+1}\b^{-\frac{2Nm}{N+1}+\frac{N}{2}}\Bigg\{N\zeta(2m)\left(\gamma+\log\left( \frac{\b}{2\pi}2^{1/N}\right)\right)-\zeta'(2m)\Bigg\}\nonumber\\
&\qquad \qquad \qquad \quad+\sum_{j=0}^{m-1}\frac{B_{2j}~\zeta(2Nm+1-2Nj)}{(2j)!}2^{N(2j-1)}\a^{2j-\frac{2Nm}{N+1}-\frac{1}{2}}.
\end{align}

To study the behavior of the series on the left-hand side of \eqref{lambert} as $\a \to 0$, we examine the series on the right-hand side of \eqref{lambert} as $\b \to \i$ using the relation $\a\b^N=\pi^{N+1}.$ 

We now invoke \eqref{f+finf} with $k=2m$, $N$ replaced by $1/N$ and $x=\b 2^{1/N}e^{i\pi j/N}$ in the expression below, then replace $j$ by $j/2$ to see that as $\b \to \infty$,
\begin{align}
&\sum_{j=-\frac{(N-1)}{2}}^{\frac{(N-1)}{2}}\left(\mathscr{F}_{2m,\frac{1}{N}}\left( \frac{i \b}{2\pi}2^{1/N}e^{\frac{i\pi j}{N}}\right)+\mathscr{F}_{2m,\frac{1}{N}}\left( -\frac{i \b}{2\pi}2^{1/N}e^{\frac{i\pi j}{N}}\right)\right)\nonumber\\
&\sim\sum_{n=1}^{\infty}(-1)^{n+1}\frac{B_{2n}}{n}\zeta\left(2m+\frac{2n}{N}\right)\left( \frac{\b}{2\pi}2^{1/N}\right)^{-2n}\sum_{j=-(N-1)}^{N-1}{\vphantom{\sum}}''e^{-\frac{i\pi nj}{N}}.\nonumber
\end{align}
From \eqref{Chebi}, for $\ell \in \mathbb{N}$,
\begin{align*}
\sum_{j=-(N-1)}^{(N-1)}{\vphantom{\sum}}''e^{-\frac{i\pi nj}{N}}=\left\{
	\begin{array}{ll}
		N(-1)^{\ell(N-1)} & \mbox{if } n=\ell N \\
		0 & \mbox{if } n\neq \ell N .
	\end{array}
\right.
\end{align*}
Hence, 
\begin{align}\label{midcal}
&\sum_{j=-\frac{(N-1)}{2}}^{\frac{(N-1)}{2}}\left(\mathscr{F}_{2m,\frac{1}{N}}\left( \frac{i \b}{2\pi}2^{1/N}e^{\frac{i\pi j}{N}}\right)+\mathscr{F}_{2m,\frac{1}{N}}\left( -\frac{i \b}{2\pi}2^{1/N}e^{\frac{i\pi j}{N}}\right)\right)\nonumber\\
&\sim N\sum_{\ell=1}^{\infty}(-1)^{\ell+1}\frac{B_{2\ell N}}{\ell N}\zeta\left(2m+2\ell\right)\left( \frac{\b}{2\pi}2^{1/N}\right)^{-2\ell N}.
\end{align}
Thus from \eqref{lambert}, \eqref{Dm} and \eqref{midcal} and the relation $\b =\frac{\pi^{(N+1)/N}}{\a^{1/N}}$, we arrive at \eqref{lambert asymptotic eqn1} upon simplification.

To obtain \eqref{lambert asymptotic eqn2}, simply put $N=1$ in \eqref{lambert asymptotic eqn1} and simplify.
\end{proof}

\section{Concluding Remarks}\label{cr}
The primary goal of this paper was to obtain functional equations for a new generalization of the Herglotz function, namely $\mathscr{F}_{k,N}(x)$. We obtained two different kinds of functional equations relating $\mathscr{F}_{k,N}(x)$ with $\mathscr{F}_{\frac{N+k-1}{N},\frac{1}{N}}(\xi/x)$, where $\xi$ is some root of unity, the first of which (Theorem \ref{hhfe}) reduces to Zagier's \eqref{fe2} when $k=N=1$.\\

(1) Although we were unable to get a three-term functional equation for $\mathscr{F}_{k,N}(x)$ similar to \eqref{fe2} and \eqref{vz2f}, we have an idea that might help suggest the form of such an equation, if it exists. This is explained below.

Note that when we let $x=1$ in \eqref{fe1} and \eqref{zaginteval}, the corresponding Herglotz functions have the same arguments \footnote{In fact, when we let $x=1$ in equations \eqref{fe1} and \eqref{zaginteval}, they exactly match since $J(1)=\frac{1}{2}\log^{2}(2).$}. Provided this phenomenon persists when we transition from the Herglotz function to the extended higher Herglotz function $\mathscr{F}_{k,N}(x)$, 
the three-term functional relation that is sought might be involving
\begin{equation*}
\mathscr{F}_{k,N}(x^{N}),\hspace{3mm} \mathscr{F}_{k,N}((1+x)^{N}),\hspace{3mm}\text{and}\hspace{3mm} \mathscr{F}_{k,N}\left(\frac{x^{N}}{(1+x)^{N}}\right),
\end{equation*}
for, not only do they reduce to $F(x)$, $F(x+1)$ and $F\left(\frac{x}{x+1}\right)$ respectively when $k=N=1$ which are the same as the ones occurring in \eqref{fe2} but they also reduce to $\mathscr{F}_{k,N}(1)$, $\mathscr{F}_{k,N}(2^{N})$ and $\mathscr{F}_{k,N}(2^{-N})$ when we let $x=1$ which are indeed the same as those occurring in \eqref{zagintevalgeneqn} for $x=1$.

Such a three-term functional relation, if/when obtained, would be a first-of-its-kind result since nowhere in the literature has there been a relation which involves like powers (which are greater than $1$) of $x$, $x+1$ and $x/(x+1)$ in the arguments of the associated functions. \\

(2) Equation \eqref{zetageneqna} involves the series, namely,
\begin{equation*}
\sum_{n=1}^{\infty}\frac{\sin(2\pi na)}{n^{2m+1}}\left(\psi\left(\tfrac{i\beta}{2\pi}(2n)^{\frac{1}{N}} e^{\frac{i\pi j}{N}}\right)+\psi\left(\tfrac{-i\beta}{2\pi}(2n)^{\frac{1}{N}} e^{\frac{i\pi j}{N}}\right)\right),
\end{equation*}
which, if we differentiate with respect to $a$ and then let $a=1$, leads to the following combination of the extended higher Herglotz functions
\begin{equation*}
\sum_{n=1}^{\infty}\frac{1}{n^{2m}}\left(\psi\left(\tfrac{i\beta}{2\pi}(2n)^{\frac{1}{N}} e^{\frac{i\pi j}{N}}\right)+\psi\left(\tfrac{-i\beta}{2\pi}(2n)^{\frac{1}{N}} e^{\frac{i\pi j}{N}}\right)\right).
\end{equation*}
This was, in fact, our motivation to study $\mathscr{F}_{k,N}(x)$. 

This suggests a further question - does there exist a transformation for a more general series
\begin{equation*}
\sum_{n=1}^{\infty}\frac{\cos(2\pi na)}{n^{k}}\left(\psi\left(\frac{i\a}{2\pi}(2n)^{N}\right)+\psi\left(\frac{-i\a}{2\pi}(2n)^{N}\right)\right),
\end{equation*}
where $k\geq1$ and $0<a\leq1$. This would then generalize Theorems \ref{genherglotzNodd} and \ref{thmtriplepole} and the general result, if obtained, would be analogous to \eqref{zetageneqna} which is the corresponding result in the setting of generalized Lambert series. This would then complete the analogy between extended higher Herglotz functions and generalized Lambert series as specified in Remark \ref{ehhfgls1}.\\

(3) As explained in the introduction, the closed-form evaluation of the integral $J(x)$ in \eqref{jx} is of interest not only from the point of view of analytic number theory but also of algebraic number theory. In the same vein, it would be worthwhile to study for which values of $k, N$ and $x$ can the integral in \eqref{defnJkN} be evaluated in closed-form. In Corollary \ref{N=1JK}, we have evaluated it in an ``almost'' closed-form since it has the constant $F_k(2)$ in its evaluation. \\

(4) Page 220 of Ramanujan's Lost Notebook \cite{lnb} contains a beautiful modular relation for $\a,\b$ positive such that $\a\b=1$, namely, if  
\begin{align*}
\phi(x):=\psi(x)+\frac{1}{2x}-\log x,
\end{align*}
and $\xi(s)$ and $\Xi(t)$ denote Riemann's functions respectively defined by \cite[p.~16, Equations (2.1.12), (2.1.14)]{titch}
\begin{align}
\xi(s)&:=(s-1)\pi^{-\tfrac{1}{2}s}\Gamma(1+\tfrac{1}{2}s)\zeta(s),\nonumber\\
\Xi(t)&:=\xi(\tfrac{1}{2}+it),\nonumber
\end{align}
then
\begin{align}
\sqrt{\a}&\left\{\frac{\gamma - \log(2 \pi \a)}{2\a} +\sum_{n=1}^{\infty}\phi(n\a)\right\}=\sqrt{\b}\left\{\frac{\gamma - \log(2 \pi \b)}{2\b} +\sum_{n=1}^{\infty}\phi(n\a)\right\}\nonumber\\
&=-\frac{1}{\pi^{3/2}}\int_{0}^{\infty}\left|\Xi\left( \frac{t}{2}\right) \Gamma\left(\frac{-1+it}{4} \right)\right|^2\frac{\cos\left(\frac{t}{2}\log\a \right)}{1+t^2}dt.\nonumber
\end{align}
The first equality in this formula seems to be in the realms of the theory of Herglotz function. A proof of it can be found, for example, in \cite{bcbad}. It would be certainly of merit to see if another proof of it could be obtained through the theory of Herglotz function.\\

(5) While the finite sum in Ramanujan's formula for $\zeta(2m+1)$, that is, \eqref{zetaodd}, involves only even-indexed Bernoulli numbers, or, in other words, even zeta values, that in \eqref{companioneqn} involves a product  of an even zeta value and an odd zeta value. 
The remaining case of having products of only odd zeta values in the summand of the finite sum is covered by our Corollary \ref{trans2m+1}.

In view of the fact that the finite sum in \eqref{zetaodd} has given rise to a nice theory of Ramanujan polynomials \cite{gmr}, \cite{msw}, it would be interesting to study the analogous polynomials stemming from \eqref{trans2m+1eqn} and \eqref{companioneqn}.\\

(6) Our method for proving Theorem \ref{hhfe} limits $k$ to be between $1$ and $N$ ($N$ inclusive). Does there exist a functional equation for $F_{k, N}(x)$ where $k$ which are greater than $N$?


\begin{center}
Acknowledgements
\end{center}
The authors sincerely thank Professor Kenneth S.~Williams for kindly sending them a copy of \cite{herglotz}. The first author's research was partially supported by SERB-DST CRG grant CRG/2020/002367. The third author is a postdoctoral fellow funded in part by the same grant. Both sincerely thank SERB-DST for the support. The third author also thanks IIT Gandhinagar for financial support.

\end{document}